\documentstyle[amscd, amstex]{amsart}

\newtheorem{thm}{Theorem}[section]
\newtheorem{lem}[thm]{Lemma}
\newtheorem{cor}[thm]{Corollary}
\newtheorem{pr}[thm]{Proposition}

\theoremstyle{definition}
\newtheorem{rem}[thm]{Remark}
\newtheorem{ex}[thm]{Example}
\newtheorem{defn}[thm]{Definition}

\newcommand{\psx}{{\bf P}_{\Sigma_X}}
\newcommand{\tta}{{\bf T}_\tau}
\newcommand{\ts}{{\bf T}_\sigma}
\newcommand{\mult}{{\rm mult}}
\newcommand{\key}{\bibitem}
\newcommand{\psd}{{\bf P}_{\Sigma_D}}
\newcommand{\ps}{{\bf P}_{\Sigma}}
\newcommand{\vol}{{\rm vol}}
\newcommand{\gr}{{\rm Gr}}
\newcommand{\res}{{\rm Res}}
\newcommand\hidot{{\raise1pt\hbox{$\scriptscriptstyle\bullet$}}}
\newcommand\lodot{{\raise.3pt\hbox{$\scriptscriptstyle\bullet$}}}
\newcommand{\dd}{{\rm d}}
\newcommand{\pp}{{\bf P}}
\newcommand{\ttt}{{\bf T}}
\newcommand{\tr}{{\rm Tr}}

\begin{document}
\title 
{Semiample hypersurfaces in toric varieties}
\author{Anvar R. Mavlyutov}
\address {Department of Mathematics \& Statistics,
University of Massachusetts, Amherst, MA, 01003, USA}
\email{anvar@@math.umass.edu}
 
\thanks{Author's work   partially supported by I. Mirkovic from  National
Science Foundation  grant  number DMS-9622863.}
\thanks{Published in Duke. Math. J. {\bf 101} (2000), 85--116.}

\begin{abstract} We study the geometry and cohomology of semiample hypersurfaces 
in toric varieties.
Such hypersurfaces generalize the  MPCP-desingulari\-zations of Calabi-Yau
ample hypersurfaces in the Batyrev mirror construction. 
We study the topological cup product on the middle cohomology of semiample hypersurfaces.
In particular, we  obtain a complete algebraic description of
the middle cohomology of regular semiample hypersurfaces in 4-dimensional simplicial toric varieties
what would be interesting for physics.  
\end{abstract}

\maketitle

\section*{Introduction}

While the geometry and cohomology of ample hypersurfaces in toric varieties have been studied  \cite{bc}, 
not much attention has been paid to semiample (i.e., ``big'' and ``nef'')
hypersurfaces defined by sections of line bundles generated by global sections with a positive 
self-intersection number. It turns out that mirror symmetric hypersurfaces in the Batyrev mirror
construction \cite{b2} are semiample, but often not ample. 
In this paper we will study semiample hypersurfaces. Such hypersurfaces
bring a geometric construction which generalizes the way of construction in \cite{b2}.
          
The purpose of this paper is to present some approaches to studying the cohomology ring 
of semiample hypersurfaces in complete simplicial toric varieties.
In particular, we  explicitly describe the ring structure on
the middle cohomology of regular semiample hypersurfaces, when the dimension of the ambient 
space is 4. 
Let us explain the main ideas of computing the topological cup product. 
The first step is to naturally 
relate the middle cohomology of the hypersurfaces to some graded ring; in our 
situation this will be done using a Gysin spectral sequence. The origins
of this ideas are in \cite{cg}, \cite{bc}.  
The second step is to use the multiplicative structure on the graded ring in
order to compute the topological cup product on the middle cohomology.
 We remark that the cup product was computed on 
the middle cohomology of smooth hypersurfaces in a projective space \cite{cg},
and this paper will generalize some of the results in \cite{cg}.

The following is a brief summary of  the paper:

In section 1 we establish notation and then introduce a geometric construction
associated with semiample divisors in complete toric varieties. At the end we give a criterion 
for a divisor to be ample (generated by global sections) in terms of intersection numbers.
This was known for simplicial toric varieties (the toric Nakai criterion), and we prove it 
for arbitrary complete toric varieties.

Section 2 studies regular semiample hypersurfaces and describes a nice stratification of such 
hypersurfaces. These hypersurfaces generalize those in the Batyrev construction \cite{b2}.

Section 3 generalizes the results of \cite {cg} on an algebraic cup product formula for residues 
of rational differential forms (from here on, the toric variety is usually simplicial). 
 It shows that there is a natural map from a graded ring (the Jacobian ring $R(f)$ \cite{bc}) to the 
middle cohomology of a quasismooth hypersurface such that the multiplicative structure on the 
ring is compatible with the topological cup product.

In section 4 we partially describe the middle cohomology  of a regular semiample hypersurface $X$;
in particular, we show that some graded pieces of the ring $R_1(f)$, considered in \cite{bc},
are imbedded into the middle cohomology of $X$. We explicitly compute the cup product on the part
coming from the ring. We should point out that, when $X$ is ample, the graded pieces of $R_1(f)$ fill
up the middle cohomology of the hypersurface, but not so in the semiample case.
 
Section 5  computes the middle 
cohomology and the cup product on it for regular semiample  hypersurfaces in a 4-dimensional 
toric variety. This is the most interesting case for physics. 
We describe the whole middle cohomology in algebraic terms, even though  $R_1(f)$ might fill
up only part of the middle cohomology. In fact, the complement to the $R_1(f)$ part is  a direct sum
of the middle cohomologies of regular ample hypersurfaces in 2-dimensional toric varieties. Hence, 
this part can also be described in terms of rings similar to $R_1(f)$. 

In section 6 we compute the Hodge numbers $h^{p,2}$ of a regular semiample 
hypersurface, and then apply the obtained formulas to the hypersurfaces in the
Batyrev mirror construction \cite{b2} to verify that, in general, the duality predicted by physicists 
does not occur for the Hodge numbers of such hypersurfaces.

Basic references on the theory
of toric varieties are \cite{f1}, \cite{o}, \cite{d}, \cite{c2}.

{\it Acknowledgments.} I would like to thank David Cox for his advice and useful
comments. I am  grateful to David Cox and David Morrison for allowing me
to use their unpublished notes for Theorems \ref{t:res} and \ref{t:cup}. I also thank the referee
for pointing out that our notion ``semiample'' is a little bit different from the common one (see Remark 
\ref{r:note} bellow).

\section{Semiample divisors}

In this section we first establish notation,
 review some  basic facts from the 
toric geometry, and then discuss a geometric construction associated with
semiample divisors in complete toric varieties. At the end of this section we will
prove  a generalization of the toric Nakai criterion for arbitrary complete toric varieties. 
As a consequence we will obtain a criterion for semiample divisors 
in terms of intersection numbers. In notation we follow \cite{bc}, \cite{c2}.
 
Let $M$ be a lattice of rank $d$, $N=\text{Hom}(M,{\Bbb Z})$ the 
dual lattice; $M_{\Bbb R}$ (resp.~$N_{\Bbb R})$ denotes
the $\Bbb R$-scalar extension of $M$ (resp.~of~$N$).
The symbol ${\bf P}_{\Sigma}$ stands for a complete  
toric variety associated with   
a finite complete  fan $\Sigma$ 
in $N_{\Bbb R}$. 
Denote by $\Sigma(k)$ the set of all $k$-dimensional cones in $\Sigma$; in particular, 
$\Sigma(1)=\{\rho_1,\dots,\rho_n\}$ is the set of $1$-dimensional
cones in $\Sigma$ with the minimal integral generators 
$e_1,\dots,e_n$, respectively.
Each 1-dimensional cone $\rho_i$ corresponds to
a torus-invariant divisor $D_i$ in ${\bf P}_\Sigma$. 
A torus-invariant Weil divisor $D=\sum_{i=1}^{n}a_iD_i$ determines
a convex polyhedron
$$\Delta_D=\{m\in M_{\Bbb R}:\langle m,e_i\rangle\geq-a_i
\text{ for all } i\}\subset M_{\Bbb R}.$$
When $D=\sum_{i=1}^{n}a_iD_i$ is Cartier, 
there is a support function 
$\psi_D : N_{\Bbb R}\rightarrow\Bbb R$ that is linear on each
cone $\sigma\in\Sigma$ and determined by some $m_\sigma\in M$:
$$\psi_D(e_i)=\langle m_\sigma, e_i\rangle=-a_i\text{ for all }
e_i\in\sigma.$$    
Since ${\bf P}_\Sigma$ is complete, a general fact is that a Cartier divisor $D$ 
(i.e., the corresponding line bundle ${\cal O}_{\ps}(D)$) is generated 
by global sections (resp. ample) if and only if $\psi_D$ is  convex (resp. strictly 
convex). 

A Cartier divisor $D$ on $\ps$ is called {\it semiample} if $D$ is 
generated by global sections and the intersection number 
$(D^d)>0$.  In complete toric varieties all ample divisors are semiample.
{}From \cite[sect.~5.3]{f1} it follows that $(D^d)=d! \vol_d(\Delta_D)$
where $\vol_d$ is the $d$-dimensional volume normalized with respect to the lattice $M$. 
So, the semiample torus-invariant divisors in complete toric varieties can be 
characterized by the two conditions that the support function
$\psi_D$ is convex and the polyhedron $\Delta_D$ has maximal dimension $d$.

\begin{rem} \label{r:note}
We should mention here that our notion ``semiample'' is a little bit different from the common one.
In \cite{ev} it is not assumed that semiample sheaves ${\cal L}$  have the additional property ${\cal L}^d>0$ 
(the Iitaka dimension is maximal).  
The author believes that the results in this section can be easily generalized for all Cartier divisors 
generated by global
sections. However, for the purpose of studying mirror symmetric hypersurfaces \cite{bc} we simply assume that 
semiample sheaves have the additional property. The same definition was used in the recent book \cite{ck}.
\end{rem}

Let us show how to construct a semiample (but not ample) divisor from an ample one.
Consider a proper birational morphism 
$\pi:{\bf P}_{\Sigma_1}\rightarrow{\bf P}_{\Sigma_2}$ 
between two complete toric varieties corresponding to a subdivision $\Sigma_1$ 
of a fan $\Sigma_2$ with an ample torus-invariant divisor $Y$ on ${\bf P}_{\Sigma_2}$.
Then the pull-back $\pi^*(Y)$ is a torus-invariant Cartier divisor with 
the same support function as the one for $Y$. Hence, $\pi^*(Y)$ is semiample and 
it is not ample if $\Sigma_1$ is different from $\Sigma_2$.    

We now show that all semiample divisors arise uniquely this way,
constructing  a complete fan $\Sigma_D$ for a semiample 
Cartier divisor $D=\sum_{i=1}^{n}a_iD_i$ using our fan 
$\Sigma$ and the convex support function $\psi_D$.
The value of the support function $\psi_D$ on 
each $d$-dimensional cone $\sigma\in\Sigma$ is determined by a unique $m_\sigma\in M$. 
We glue together those maximal dimensional cones in $\Sigma$ that have the same 
$m_\sigma$. The glued set is again a convex rational polyhedral cone, 
and one can show that
this cone is strongly convex using the fact that $\Delta_D$ has maximal dimension $d$.
The set of these strongly convex cones with its faces comprise a new complete
fan $\Sigma_D$ in $N_{\Bbb R}$. 
This construction is independent of the equivalence relation on the divisors:
if we change the divisor $D$ to  a linearly equivalent one, the fan $\Sigma_D$
will remain the same. 

The fan $\Sigma_D$ is exactly the normal fan of $\Delta_D$. Indeed, by construction, 
$\psi_D$ is strictly convex with respect to $\Sigma_D$.
On the other hand, since $D$ is generated  by global sections, the support function 
$\psi_D$ coincides with the function of $\Delta_D$ \cite[sect.~3.4]{f1}:
$$\psi_D(n)=\min_{m\in\Delta_D}\langle m,n\rangle.$$
Theorem 2.22 \cite{o} implies that $\Sigma_D$ is the normal fan of $\Delta_D$.

Notice that $\Sigma$ is a refinement of $\Sigma_D$. So,  the sets of 
1-dimensional cones of the fans are related by $\Sigma_D(1)\subset\Sigma(1)$, and we have
a proper birational morphism 
$\pi:{\bf P}_\Sigma\rightarrow{\bf P}_{\Sigma_D}$ between the two toric varieties.
Any proper morphism determines a push-forward map 
$\pi_*:A_{d-1}({\bf P}_\Sigma)\rightarrow A_{d-1}({\bf P}_{\Sigma_D})$ on 
the Chow group, that takes the class of an irreducible divisor $V$ to
the class  $\deg(V/\pi(V))[\pi(V)]$ if $\pi(V)$ has the same dimension as $V$
and to 0 otherwise. Now apply the push-forward map to our semiample divisor:      
$$\pi_*[D]=\sum\limits_{i=1}^{n}a_i\pi_*[D_i]=[\sum_{\rho_i\in\Sigma_D(1)}a_i\pi(D_i)],$$ 
because
$D_i$ maps birationally onto its image when $\rho_i\in\Sigma_D$, and $\dim\pi(D_i)<\dim D_i$
in all other cases.
The divisors $\pi(D_i)$ for $\rho_i\in\Sigma_D(1)$ are torus-invariant corresponding to
the 1-dimensional cones in $\Sigma_D$. The support function of
the Weil divisor $\pi_*(D):=\sum_{\rho_i\in\Sigma_D(1)}a_i\pi(D_i)$
 coincides with $\psi_D$ which 
is strictly convex with respect to the fan $\Sigma_D$. Hence, the divisor class $\pi_*[D]$
 is ample. 
For the birational map 
$\pi:{\bf P}_\Sigma\rightarrow{\bf P}_{\Sigma_D}$
we also have a commutative diagram \cite[sec. 3.4]{f1}, \cite{f2}:
\[ \begin{CD}
A_{d-1}({\bf P}_\Sigma)@>\pi_*>>A_{d-1}({\bf P}_{\Sigma_D})\\
@AAA  @AAA \\
\text{Pic}({\bf P}_\Sigma) @<<\pi^*< \text{Pic}({\bf P}_{\Sigma_D}),
\end{CD} \]
where the vertical maps are inclusions.
Since the support functions for the Cartier divisors $D$ and $\pi_*(D)$ coincide, 
the pull-back
$\pi^*\pi_*[D]$ is exactly the divisor class $[D]$. Thus, we have the following useful  result.
          
\begin{pr}\label{p:sem} 
Let $\bold P_\Sigma$ be a complete toric variety with a semiample divisor class $[D]\in A_{d-1}(\ps)$.
There exists a unique complete toric variety $\psd$ with a toric birational map
$\pi:{\bf P}_\Sigma\rightarrow{\bf P}_{\Sigma_D}$, such that   $\Sigma$ is a subdivision of 
$\Sigma_D$, $\pi_*[D]$ is ample and $\pi^*\pi_*[D]=[D]$. Moreover, if 
$D=\sum_{i=1}^{n}a_i D_i$ is torus-invariant, then $\Sigma_D$ is the normal fan of $\Delta_D$.
\end{pr}
 
\begin{rem}\label{r:cor} Since  the fan $\Sigma_D$ is  the normal fan of $\Delta_D$,
there is a one-to-one correspondence between the $k$-dimensional cones of $\Sigma_D$ and
$(d-k)$-dimensional faces of $\Delta_D$. Note, however, that while $\Sigma_D$ is canonical with 
respect to the equivalence relation on the divisors, the polyhedron 
$\Delta_D$ is only canonical up to translation. 
\end{rem}

We next study the intersection theory for the semiample divisors in 
complete toric varieties. 
Any toric variety $\ps$ is a disjoint union of its orbits by the action of the 
torus 
${\bf T}=N\otimes{\Bbb C}^*$ that sits naturally inside 
$\ps$. Each orbit ${\bf T}_\sigma$ is a torus corresponding to a cone 
$\sigma\in\Sigma$. The closure of each orbit ${\bf T}_\sigma$ is
again a toric variety denoted $V(\sigma)$.

\begin{lem}\label{l:int}
 If $D$ is a semiample divisor on a complete toric variety  
$\ps$, then the intersection number
$(D^k\cdot V(\sigma))>0$ for any $\sigma\in\Sigma(d-k)$
contained in a cone of $\Sigma_D(d-k)$, 
and $(D^k\cdot V(\sigma))=0$ for all other $\sigma\in\Sigma(d-k)$.
\end{lem}

\begin{pf}
We can assume that $D=\sum_{i=1}^{n}a_i D_i$, which gives a support function $\psi_D$
determined on each cone by some $m_\sigma$: $\psi_D(n)=\langle m_\sigma,n\rangle$ for all
$n\in\sigma$.
Since $D$ is generated by global sections, for a fixed
$\sigma\in\Sigma(d-k)$ we have \cite[sec. 5.3]{f1}:
\begin{equation}
\vol_k(\Delta_D\cap(\sigma^\perp+m_\sigma))=(\frac{D^k}{k!}\cdot V(\sigma)).
\label{e:int}
\end{equation}
By Remark \ref{r:cor},  there is a one-to-one correspondence between the  cones of $\Sigma_D$ and the
faces of $\Delta_D$. 
Let $\sigma_\Gamma$ be the minimal
cone in $\Sigma_D$, corresponding to a face $\Gamma$ of $\Delta_D$  and containing $\sigma$. We claim that
$$\Gamma=\Delta_D\cap(\sigma^\perp+m_\sigma).$$
Indeed, since $\psi_D$ is strictly convex with respect to $\Sigma_D$, from  Lemma 2.12 
\cite{o} we have  
$$\Gamma=\{m\in\Delta_D:\langle m,n\rangle=\psi_D(n)\mbox{ for all }n\in\sigma_\Gamma\},$$  
whence $m\in\Gamma$ implies $\langle m-m_\sigma,n\rangle=0$ for all $n\in\sigma$.
Conversely, suppose $m\in\Delta_D$ and $(m-m_\sigma)\in\sigma^\perp$. The first condition 
implies
$\langle m,n\rangle\geq\psi_D(n)$ for all $n$ from the strongly convex cone $\sigma_\Gamma$,
while the second one gives a point in the interior of $\sigma_\Gamma$  (by the minimal 
choice of 
this cone) for which $m$ and $\psi_D$ have the same values. Hence, $m$ and $\psi_D$ have 
the same 
values on $\sigma_\Gamma$, and the claimed equality of the sets follows.

Now, the lemma follows from the fact that $\vol_k(\Gamma)>0$ if and only if 
$\dim\sigma_\Gamma=d-k$.
\end{pf} 

\begin{rem}\label{r:int}
The above lemma  provides another way of constructing the fan $\Sigma_D$,
by gluing the $d$-dimensional cones in $\Sigma$ along those facets $\tau$ for which
$(D\cdot V(\tau))=0$.
\end{rem}
           
We will now give necessary and sufficient conditions for a Cartier divisor on a complete 
toric variety 
to be ample, generated by global sections or semiample. This is a generalization of the 
toric Nakai criterion proved for nonsingular toric varieties in Theorem~2.18 \cite{o}.

\begin{thm}\label{t:cr} 
Let $\ps$ be a $d$-dimensional complete toric variety and let $D$ be a Cartier divisor on 
$\ps$. 
Then

{\rm(i)} $D$ is generated by global sections if and only if $(D\cdot V(\tau))\ge0$ for any 
$\tau\in\Sigma(d-1)$.

{\rm(ii)} $D$ is ample if and only if $(D\cdot V(\tau))>0$ for any $\tau\in\Sigma(d-1)$.
\end{thm}

\begin{pf} Without loss of generality we can assume that $D$ is torus-invariant.

(i) If $D$ is generated by global sections, then the required condition follows from equation 
(\ref{e:int}).
Conversely, the torus-invariant divisor $D$ has the support function $\psi_D$, and it suffices 
to show 
that $\psi_D$ is convex.
Here, we use a trick.
Consider a nonsingular subdivision $\Sigma'$ of the fan $\Sigma$ 
and the corresponding toric morphism $f:{\bf P}_{\Sigma'}\rightarrow\ps$.
Then the support function of the pull-back divisor $f^*(D)$ coincides with $\psi_D$.
So, we just need to show that $f^*(D)$ is generated by global sections.
By Example 2.4.3 \cite{f2}, we have $(f^*(D)\cdot V(\tau'))=(D\cdot f_*(V(\tau')))$, 
where $V(\tau')$ is the closure of the 1-dimensional orbit corresponding to $\tau'\in\Sigma'(d-1)$. 
If the smallest cone in $\Sigma$ containing $\tau'$ is $d$-dimensional, then 
the image of $V(\tau')$ is a point, implying that the above intersection number vanishes.
Otherwise, $\tau'$ is contained in some $\tau\in\Sigma(d-1)$, in which case 
$f_*(V(\tau'))=V(\tau)$. 
So, in either case, by the given condition in (i), the intersection number  $(D\cdot f_*(V(\tau')))$ is 
nonnegative.
Following the proof of Theorem 2.18 \cite{o} we get that $(f^*(D)\cdot V(\tau'))\ge0$ for any 
$\tau'\in\Sigma'(d-1)$ implies $f^*(D)$ is generated by global sections. 

(ii) If $D$ is ample, then the required condition follows from Lemma \ref{l:int} or, 
more generally, 
from the Nakai criterion for arbitrary complete varieties \cite[chap.~I,~Theorem~5.1]{h}, \cite{k}.  

Conversely, by part (i), the divisor $D$ is generated by global sections. 
We will show that $D$ is semiample and the fan $\Sigma$ is exactly the fan $\Sigma_D$ associated 
with the semiample 
divisor. Then, by Proposition \ref{p:sem}, the desired result will follow.
{}From equation (\ref{e:int}) and the given condition 
it follows that the polyhedron $\Delta_D$ intersects different lines, corresponding 
to $\tau\in\Sigma(d-1)$, in more than one point. These lines can not lie in a hyperplane of 
$M_{\Bbb R}$, 
because $\Sigma$ is complete. Therefore, $\Delta_D$ is maximal dimensional, implying that $D$ is 
semiample.
By the Remark \ref{r:int} and the given condition, the fan $\Sigma$ coincides with $\Sigma_D$.
Thus, Proposition \ref{p:sem} implies that $D$ is ample.
\end{pf}

\begin{cor}
Let $\ps$ be a complete  toric variety. Then a Cartier divisor $D$ on $\ps$
 is semiample if and only if $(D^d)>0$ and $(D\cdot V(\tau))\geq0$ for any 
$\tau\in\Sigma(d-1)$.
\end{cor}

\begin{rem} 
In  Mori's theory, Theorem \ref{t:cr}(ii) above and Proposition (1.6) of \cite{r}  
imply that $D$ is ample if and only if 
$(D\cdot({\Bbb NE}(\ps)\setminus\{0\}))>0$, 
where ${\Bbb NE}(\ps)$ is the cone coming from effective 1-cycles.
Also, by part (i) of Theorem \ref{t:cr}, the pseudo-ample cone ${\Bbb PA}(\ps)$ is spanned by the 
divisors generated by global sections.
For details see \cite{r}, \cite[sect.~2.5]{o}. 
\end{rem}

\section{Regular semiample hypersurfaces}

Next we shall apply results from the previous section to describe a stratification of 
regular  semiample hypersurfaces in a complete toric variety $\ps$. 
The following definition has appeared in \cite{b2}.

\begin{defn}
A hypersurface $X\subset\ps$ is called {\it $\Sigma$-regular} if $X\cap\ts$ is empty or 
a smooth subvariety of codimension 1 in $\ts$ for any $\sigma\in\Sigma$.
\end{defn}

\begin{rem}
Proposition 6.8 \cite{d} says that a hypersurface $X\subset\ps$ defined by a general section of 
a line bundle generated by global sections is $\Sigma$-regular.
When it is clear from the context, we simply say that a hypersurface is regular.
\end{rem}

\begin{lem}\label{l:ir} Let $X$ be a semiample hypersurface in a complete toric variety $\ps$, such that
$\dim\ps\ge2$. Then

{\rm(i)} $X$  is connected, and

{\rm(ii)} $X$  is irreducible if $X$ is $\Sigma$-regular. 
\end{lem}

\begin{pf} {\rm(i)}
Consider an effective torus-invariant divisor $D$  equivalent to the divisor $X$.
Since ${\cal O}_{\ps}(D)$ is generated by global sections, choosing a basis of the space 
$H^0(\ps,{\cal O}_{\ps}(D))$  gives a mapping  
$$\varphi_D:\ps\rightarrow{\Bbb P}^{r-1},$$
where $r=h^0(\ps,{\cal O}_{\ps}(D))=\mbox{Card}(\Delta_D\cap M)$.
By Exercise on p.~73 \cite[sect.~3.4]{f1}, the image of $\varphi_D$ has dimension equal to 
$\dim\Delta_D$. Since $D$ is semiample, we get that $\dim\Delta_D=\dim\ps\ge2$.
{}From Theorem 2.1 \cite{fl} it follows that 
every divisor in the linear system $|D|$ is connected. In particular, $X$ is connected.

{\rm(ii)}
To prove that $X$ is irreducible we argue as follows. Consider a nonsingular subdivision 
$\Sigma'$
of the fan $\Sigma$ and the corresponding morphism $p:{\bf P}_{\Sigma'}\rightarrow\ps$.
It follows from Proposition 3.2.1 \cite{b2} that  
$p^{-1}(X)$ is  a $\Sigma'$-regular hypersurface that supports a semiample divisor  $p^*(X)$.
By the previous part, $p^{-1}(X)$ is a smooth connected hypersurface which
 must be irreducible. Therefore, $X$ is irreducible. 
\end{pf}

\begin{pr}\label{p:reg}
 If $X\subset\ps$ is a $\Sigma$-regular semiample  hypersurface 
with the associated morphism $\pi:\ps\rightarrow\psx$ for the divisor class $[X]\in A_{d-1}(\ps)$ 
from Proposition \ref{p:sem},
then  $Y:=\pi(X)$ is a $\Sigma_X$-regular ample hypersurface, and $X=\pi^{-1}(Y)$. 
\end{pr}

\begin{pf}   
{}From Lemma \ref{l:ir}(ii) we know that $X$ is irreducible. Since $X$ is $\Sigma$-regular,
it maps birationally onto its image, implying $\pi_*[X]=[\pi(X)]$. Therefore, by Proposition 
\ref{p:sem}, the hypersurface $Y=\pi(X)$ is ample. 
  
Let us now show that $Y$ misses the 0-dimensional orbits in $\psx$.
Consider the 1-dimensional orbit closure $V(\tau_0)\subset\psx$ corresponding to a cone 
$\tau_0\in\Sigma_X(d-1)$,
and take a cone $\tau\in\Sigma(d-1)$ that lies in $\tau_0$. 
Since $X$ is $\Sigma$-regular,   
$$\mbox{Card}(X\cap\tta)=\mbox{Card}(X\cap V(\tau))=(X\cdot V(\tau)).$$ 
We also know that the orbit $\tta$ maps  onto ${\bf T}_{\tau_0}$, hence,
$$(Y\cdot V(\tau_0))\ge\mbox{Card}(Y\cap V(\tau_0))\ge\mbox{Card}(Y\cap{\bf T}_{\tau_0})
\ge\mbox{Card}(X\cap\tta)=(X\cdot V(\tau)).$$
By Example 2.4.3 \cite{f2}, we have $(Y\cdot V(\tau_0))=(X\cdot V(\tau))$, whence the above 
inequalities are equalities. Therefore, the hypersurface $Y$ intersects transversally the orbit 
${\bf T}_{\tau_0}$ and  does not intersect the points in the compliment 
$V(\tau_0)\setminus{\bf T}_{\tau_0}$, corresponding to the $d$-dimensional cones in $\Sigma_X$
that contain $\tau_0$.
Thus, we have shown that $Y$ misses
all 0-dimensional  orbits in $\psx$.

One can easily show  $X=\pi^{-1}(Y)$ from the facts that
$X$ and $Y=\pi(X)$ are irreducible, and that $Y$ misses 
the 0-dimensional orbits. 
Finally, for arbitrary $\sigma_0\in\Sigma_X$ take $\sigma\in\Sigma$, contained in $\sigma_0$,
 of the same dimension.
Then we have an isomorphism $\ts\cong{\bf T}_{\sigma_0}$ inducing another isomorphism
$$X\cap\ts=\pi^{-1}(Y)\cap\ts\cong Y\cap{\bf T}_{\sigma_0}.$$ 
So, the $\Sigma$-regularity of $X$ implies that $\pi(X)$ is $\Sigma_X$-regular.
\end{pf}

\begin{rem}\label{r:mir}
By construction in \cite{b2}, 
the MPCP-desingularizations $\widehat{Z}$ of regular projective hypersurfaces
$\overline{Z}$ in a  toric Fano variety $\pp_\Delta$ are regular semiample  hypersurfaces.  
The above proposition shows that if we start with an arbitrary regular semiample  hypersurface, then
we come up with a similar picture. 
\end{rem}

Let us note that a regular ample hypersurface $Y$ in 
a complete toric variety $\pp$ will intersect all orbits transversally, except for 0-dimensional orbits. 
Such a hypersurface is called {\it nondegenerate} in \cite{bc}, \cite{dk}. 
Also, in this case a hypersurface in the torus $\ttt$ isomorphic to 
the affine hypersurface $Y\cap{\bf T}$ in  $\ttt$
will be  called {\it nondegenerate}. Such a hypersurface satisfies the following property.

\begin{lem} \cite{dk}\label{l:lef}
Let $Z$ be a nondegenerate affine hypersurface in the torus $\ttt$, then
the natural map 
$H^i(\ttt)\rightarrow H^i(Z)$, induced by the inclusion is an isomorphism for $i<\dim\ttt-1$ and 
an injection for $i=\dim\ttt-1$. 
\end{lem}

Like in \cite{b2}, by Proposition \ref{p:reg}, we get a nice stratification
of a semiample regular hypersurface $X\subset\ps$ in terms of nondegenerate affine
hypersurfaces.
Let $Y=\pi(X)$, then  $X=\pi^{-1}(Y)$.
Using the standard description of a toric blow-up, we obtain
\begin{equation}\label{e:str}
X\cap\ts\cong (Y\cap{\bf T}_{\sigma_0})\times({\Bbb C}^*)^{\dim\sigma_0-\dim\sigma},
\end{equation}
where $\sigma_0\in\Sigma_X$ is the smallest cone containing $\sigma\in\Sigma$.

\section{A cup product formula for quasismooth hypersurfaces} 

The purpose of this section is to give a generalization of the algebraic cup product 
formula for the residues of rational forms presented in \cite{cg}. 
In this section we assume that 
$\pp$ is a complete simplicial toric variety. 
Such a  toric variety  has a homogeneous coordinate ring
$S={\Bbb C}[x_1,\dots,x_n]$ with variables $x_1,\dots,x_n$
corresponding to the irreducible torus-invariant divisors   $D_1,\dots,D_n$ \cite{c1}.
This ring is graded by the Chow group $A_{d-1}(\pp)$: $\deg(\prod_{i=1}^n x_i^{a_i})=
[\sum_{i=1}^n a_i D_i]$. 
Furthermore, if $\cal L$ is a line bundle on $\pp$, then for 
$\beta=[{\cal L}]\in A_{d-1}(\pp)$ one has an isomorphism $H^0(\pp,{\cal L})\cong S_\beta$.
So, the homogeneous polynomials in $S_\beta$ identified with the global sections of $\cal L$ 
determine  hypersurfaces in the toric variety $\pp$. 

\begin{defn}\cite{bc}
A hypersurface $X\subset\pp$ defined by a homogeneous polynomial $f\in S_\beta$ is called 
{\it quasismooth} if  $\frac{\partial f}{\partial x_i}$, $1\le i\le n$, do not
vanish simultaneously on $\pp$. 
\end{defn}

\begin{defn}\cite{bc}\label{d:om}
Fix an integer basis $m_1,\dots,m_d$ for the lattice $M$. Then given subset
$I=\{i_1,\dots,i_d\}\subset\{1,\dots,n\}$, denote 
$\det(e_I)=\det(\langle m_j,e_{i_k}\rangle_{1\le j,k\le d})$, 
$dx_I=dx_{i_1}\wedge\dots\wedge dx_{i_d}$ and $\hat{x}_I=\prod_{i\notin I}x_i$.
Define the $n$-form $\Omega$ by the formula
 $$\Omega=\sum_{|I|=d}\det(e_I)\hat{x}_I dx_I,$$
where the sum is over all $d$ element subsets $I\subset\{1,\dots,n\}.$
\end{defn}

Let $X\subset\pp$ be a quasismooth hypersurface 
defined by $f\in S_\beta$. For $A\in S_{(a+1)\beta-\beta_0}$ 
(here, $\beta_0:=\sum_{i=1}^n\deg(x_i)$) 
consider a rational $d$-form 
$$\omega_A:=A\Omega/f^{a+1}\in H^0(\pp,\Omega^d_{\pp}((a+1)X)).$$ 
This form gives a class in $H^d(\pp\setminus X)$, and by the residue map 
$$\res: H^d(\pp\setminus X)\rightarrow H^{d-1}(X)$$
we get $\res(\omega_A)\in H^{d-1}(X)$. We will need an explicit algebraic formula for
the Hodge component $\res(\omega_A)^{d-1-a,a}$ in \v{C}ech cohomology. 
 
Denote $f_i=\frac{\partial f}{\partial x_{i}}$ and let
$U_i=\{x\in\pp: f_i(x)\ne0\}$  for $i=1,\dots,n$. If $X$ is a quasismooth
hypersurface, then ${\cal U}=\{U_i\}_{i=1}^n$ is an open cover of  $\pp$.

The next two theorems with their proofs 
are  corrected and generalized versions of unpublished results of D.~Cox and D.~Morrison.

\begin{thm} \label{t:res}
Let $X\subset\pp$ be a quasismooth hypersurface 
defined by $f\in S_\beta$ and $A\in S_{(a+1)\beta-\beta_0}$, 
$\beta_0=\sum_{i=1}^n\deg(x_i)$.
Then under the natural map
$$\check{H}^a({\cal U}|_X,\Omega_X^{d-1-a})\rightarrow H^a(X,\Omega_X^{d-1-a})\cong H^{d-1-a,a}(X)$$   
the component $\res(\omega_A)^{d-1-a,a}$\! corresponds to the \v{C}ech cocycle
$c_a\fracwithdelims\{\}{AK_{i_a}\cdots K_{i_0}\Omega}{f_{i_0}\cdots f_{i_a}}
_{i_0\dots i_a}\!,$
where $c_a=\frac{1}{a!}(-1)^{d-1+a(a+1)/2}$,  and 
$K_i$  is the contraction operator $\frac{\partial}{\partial x_{i}}\lrcorner$.
\end{thm}

\begin{pf} 
The residue map can be calculated in hypercohomology using 
the commutative diagram  
$$
\begin{CD}
H^d(\pp\setminus X)@>\res>> H^{d-1}(X)\\
@AAA  @AAA \\
{\Bbb H}^d(\Omega_{\pp}^\hidot(\log X))@>\res>>{\Bbb H}^{d-1}(\Omega^{\hidot}_X),
\end{CD}
$$
where the vertical maps are isomorphisms.
As in \cite{cg} we can work in the \v{C}ech-deRham complex 
$C^\hidot({\cal U}, \Omega^\hidot(*X))$ with arbitrary algebraic singularities along $X$, where 
${\cal U}=\{U_i\}_{i=1}^n$. Then we can apply the arguments of \cite{cg} on pp.~58-62 almost 
without any change.
We only need to check that  
\begin{equation}\label{e:equ}
{\rm d}f\wedge\Omega\equiv0\mbox{ modulo multiples of } f
\end{equation} 
for part (i) of the lemma on p.~60 in \cite{cg}.
But ${\rm d}f\wedge({\rm d}x_1\wedge\dots\wedge{\rm d}x_n)=0$, and,
by Lemma 6.2 \cite{c3}, 
$\Omega=\theta_1\lrcorner\cdots\lrcorner\theta_{n-d}\lrcorner
({\rm d}x_1\wedge\dots\wedge{\rm d}x_n)$ for some Euler vector fields $\theta_i$.
The equivalence (\ref{e:equ}) can be obtained repeatedly applying the following argument.
If ${\rm d}f\wedge\omega\equiv 0$ mod $f$ for some form $\omega$, 
and $\theta$ is an  Euler vector field,
then $0\equiv\theta\lrcorner({\rm d}f\wedge\omega)=(\theta\lrcorner{\rm d}f)\omega-
{\rm d}f\wedge(\theta\lrcorner\omega)$.
Since $\theta\lrcorner{\rm d}f=\theta(f)=\theta(\beta)f$ (see the proof of Proposition 5.3
\cite{c3}), we get 
${\rm d}f\wedge(\theta\lrcorner\omega)\equiv0$. 
Thus, the lemma on p.~60 of \cite{cg} is true in our situation. 
 The rest of the arguments applies without change, and the theorem is proved.
Let us remark that the constructive proof of \cite{cg} implies that
$c_a\fracwithdelims\{\}{AK_{i_a}\cdots K_{i_0}\Omega}{f_{i_0}\cdots f_{i_a}}
_{i_0\dots i_a}$
is actually a \v{C}ech cocycle. 
\end{pf}

\begin{defn} 
For $\beta=[\sum_{i=1}^n b_i D_i]\in A_{d-1}(\pp)$ and a multi-index $I=(i_0,\dots,i_d)$ 
viewed as an ordered subset of $\{1,\dots,n\}$, we introduce
a constant $c_I^\beta$ which is the determinant of the $(d+1)\times(d+1)$ matrix 
obtained from
($\langle m_j, e_{i_k}\rangle_{1\le j\le d, i_k\in I})$ by adding the first row 
$(b_{i_0},\dots,b_{i_d})$, where  $m_1,\dots,m_d$ is the fixed integer basis
of the lattice $M$ as in Definition \ref{d:om}. 
One can easily check that $c_I^\beta$ is well defined.
\end{defn}

\begin{thm}\label{t:cup}
Let $X\subset\pp$ be a quasismooth hypersurface 
defined by $f\in S_\beta$, and  suppose $a+b=d-1$, 
$\omega_A=\frac{A\Omega}{f^{a+1}}$, $\omega_B=\frac{B\Omega}{f^{b+1}}$
for $A\in S_{(a+1)\beta-\beta_0}$, $B\in S_{(b+1)\beta-\beta_0}$.
Then under the composition 
\[\minCDarrowwidth0.7cm
\begin{CD}
H^a(X,\Omega_X^b)\otimes H^b(X,\Omega_X^a)@>\cup>> H^{d-1}(X,\Omega_X^{d-1})@>\delta>>
H^d(\pp,\Omega_{\pp}^d)
\end{CD}
\]
(here, $\delta$ is the coboundary map in the Poincar\'e residue sequence)
we have that $\delta(\res(\omega_A)^{ba}\cup\res(\omega_B)^{ab})$ 
is represented by the \v{C}ech cocycle
$$c_{ab}\fracwithdelims\{\}{A B c_I^\beta\hat{x}_I\Omega}{f_{i_0}\cdots f_{i_d}}
_I\in\check{H}^d({\cal U},\Omega_{\pp}^d),$$
where $I=(i_0,\dots,i_d)$ and $c_{ab}=\frac{(-1)^{a(a+1)/2+b(b+1)/2+a^2+d-1}}{a!b!}$.
\end{thm}

\begin{pf}
As in \cite{cg} on p.~63, by Theorem \ref{t:res}  we see that the residue product 
$\res(\omega_A)^{ba}\cup\res(\omega_B)^{ab}$ is represented by the cocycle 
$$\psi=\tilde{c}_{ab}\biggl\{\frac
{AK_{i_a}\cdots K_{i_0}\Omega}{f_{i_0}\cdots f_{i_a}}\wedge
\frac{BK_{i_{d-1}}\cdots K_{i_a}\Omega}{f_{i_a}\cdots f_{i_{d-1}}}
\biggr\}
_{i_0\dots i_{d-1}}\in C^{d-1}({\cal U}|_X,\Omega_X^{d-1}),$$ 
where $\tilde{c}_{ab}=(-1)^{a^2}c_a c_b=\frac{(-1)^{a(a+1)/2+b(b+1)/2+a^2}}{a!b!}$.
To calculate the coboundary of this cocycle we use the following commutative diagram:
$$\minCDarrowwidth0.8cm
\begin{CD} 
0 @>>> C^d({\cal U},\Omega^d_{\pp}) @>>> C^d({\cal U},\Omega_{\pp}^d(\log X)) @>\res>>
 C^{d}({\cal U}|_X,\Omega^{d-1}_X) \\
@.    @AAA                               @AAA             @AAA   \\
0 @>>> C^{d-1}({\cal U},\Omega^d_{\pp}) @>>> C^{d-1}({\cal U},\Omega_{\pp}^d(\log X)) 
@>\res>> C^{d-1}({\cal U}|_X,\Omega^{d-1}_X). 
\end{CD}
$$
Lift the cocycle $\psi$ to 
$$\tilde{\psi}=\tilde{c}_{ab}\biggl\{\frac
{AK_{i_a}\cdots K_{i_0}\Omega}{f_{i_0}\cdots f_{i_a}}\wedge
\frac{BK_{i_{d-1}}\cdots K_{i_a}\Omega}{f_{i_a}\cdots f_{i_{d-1}}}
\wedge\frac{{\rm d}f}{f}\biggr\}
_{i_0\dots i_{d-1}}\in  C^{d-1}({\cal U},\Omega_{\pp}^d(\log X)).$$
{}From the diagram we can see that changing of the numerator by a multiple of $f$ will not affect 
the image of $\psi$  in  $\check H^d({\cal U},\Omega_{\pp}^d)$.
Hence, we need to compute 
$K_{i_{a}}\cdots K_{i_0}\Omega\wedge K_{i_{d-1}}\cdots K_{i_a}\Omega\wedge{\rm d}f$ modulo multiples of $f$.  
First we will show 
\begin{equation}\label{e:mod}
K_{i_{a}}\cdots K_{i_0}\Omega\wedge K_{i_{d-1}}\cdots K_{i_a}
\Omega\wedge{\rm d}f\equiv(\mbox{some
function})\cdot\Omega\mbox{ mod } f.
\end{equation}

As in the proof of the previous theorem, we can write 
$\Omega=E\lrcorner\dd{\bf x}$, where $E$ is a wedge of some Euler vector fields and
$\dd{\bf x}=\dd x_1\wedge\dots\wedge\dd x_n$. Denote 
$\dd u=\dd x_{i_0}\wedge\dots\wedge\dd x_{i_{a-1}}$, 
$\dd v=\dd x_{i_{a+1}}\wedge\dots\wedge\dd x_{i_{d-1}}$, and
$\dd w=\wedge_{i\notin I_0}\dd x_{i}$, where $I_0=(i_0,\dots,i_{d-1})$.
Then $\dd{\bf x}=\pm\dd u\wedge\dd x_{a}\wedge\dd v\wedge\dd w$.
Now compute: 
\begin{align*}
K_{i_a}\cdots K_{i_0}\Omega&=K_{i_a}\cdots K_{i_0}(E\lrcorner\dd{\bf x})=
\pm E\lrcorner(K_{i_a}\cdots K_{i_0}\dd u\wedge\dd x_{a}\wedge\dd v\wedge\dd w) \\
&=\pm E\lrcorner(\dd v\wedge\dd w)=
\pm((E\lrcorner\dd v)\wedge\dd w+(-1)^{(d-a-1)(n-d)}\dd v(E\lrcorner\dd w)).
\end{align*}
Similarly, 
$$K_{i_{d-1}}\cdots K_{i_a}\Omega=\pm((E\lrcorner\dd u)\wedge\dd w+(-1)^{a(n-d)}
\dd u(E\lrcorner\dd w)).$$
Since $\dd w\wedge\dd w=0$, we get 
\begin{align*}
&K_{i_a}\cdots K_{i_0}\Omega\wedge K_{i_{d-1}}\cdots K_{i_a}\Omega=\pm(E\lrcorner\dd w)
\biggl((E\lrcorner\dd v)\wedge\dd w\wedge(-1)^{a(n-d)}\dd u \\
&+(-1)^{(d-a-1)(n-d)}\dd v\wedge
(E\lrcorner\dd u)\wedge\dd w+(-1)^{(d-1)(n-d)}\dd v\wedge\dd 
u\wedge(E\lrcorner\dd w)\biggr)\\
&=\pm(E\lrcorner\dd w)(E\lrcorner(\dd v\wedge\dd u\wedge\dd w))=
\pm(E\lrcorner\dd w)(E\lrcorner K_{i_a}\lrcorner\dd{\bf x}) \\ 
&=\pm(E\lrcorner\dd w)
(K_{i_a}\lrcorner(E\lrcorner\dd{\bf x}))=\pm(E\lrcorner\dd w)(K_{i_a}\Omega).
\end{align*}

{}From equation (\ref{e:equ}) we know that $\Omega\wedge\dd f\equiv0$ modulo multiples of 
$f$.
Applying the contraction operator $K_{i_a}$ to this identity we  obtain
$(K_{i_a}\Omega)\wedge\dd f\equiv\pm f_{i_a}\Omega$, whence equation (\ref{e:mod}) follows:
$$ K_{i_a}\cdots K_{i_0}\Omega\wedge K_{i_{d-1}}\cdots K_{i_a}\Omega\wedge{\rm d}f=\pm(E\lrcorner
\dd w)(K_{i_a}\Omega)\wedge\dd f\equiv\pm(E\lrcorner\dd w)f_{i_a}\Omega.$$

We next  claim that  
\begin{equation}\label{e:mod1}
K_{i_a}\cdots K_{i_0}\Omega\wedge K_{i_{d-1}}\cdots K_{i_a}\Omega\wedge{\rm d}f
\equiv(-1)^{d-1}\det(e_{I_0})\hat{x}_{I_0} f_{i_a}\Omega\mbox{ mod } f.
\end{equation}
Examine the coefficient of $\dd x_{I_0}=\dd x_{i_0}\wedge\dots\wedge\dd x_{i_{d-1}}$ in 
the left hand side.
The only place $\dd x_{I_0}$ can occur is in
\begin{align*}
&(K_{i_a}\cdots K_{i_0}\det(e_{I_0})\hat{x}_{I_0}\dd x_{I_0})\wedge
(K_{i_{d-1}}\cdots K_{i_a}\det(e_{I_0})\hat{x}_{I_0}\dd x_{I_0})\wedge f_{i_a}\dd x_{i_{a}}\\
&=\det(e_{I_0})^2\hat{x}_{I_0}^2 \dd x_{i_{a+1}}\wedge\dots\wedge\dd x_{i_{d-1}}
\wedge(-1)^{a(d-a)}
\dd x_{i_0}\wedge\dots\wedge\dd x_{i_{a-1}}\wedge f_{i_a}\dd x_{i_{a}} \\   
&=(-1)^{d-1}\det(e_{I_0})\hat{x}_{I_0}f_{i_a}(\det(e_{I_0})\hat{x}_{I_0}\dd x_{I_0}).
\end{align*}
{}From here, equation (\ref{e:mod1}) follows, 
because $\Omega=\sum_{|I|=d}\det(e_I)\hat{x}_I dx_I$, and the left hand side of (\ref{e:mod1}) is
$(\mbox{some function})\cdot\Omega$ modulo multiples of $f$. 

Returning to the calculation of the coboundary $\delta(\psi)$,
by equation (\ref{e:mod1}), we have  
$$\tilde{\psi}\equiv\frac{(-1)^{d-1}\tilde{c}_{ab}}{f}
\fracwithdelims\{\}{A B\det(e_{I_0})
\hat{x}_{I_0}\Omega}{f_{i_0}\cdots f_{i_{d-1}}}_{i_0\dots i_{d-1}},$$
so that 
$$\delta(\psi)\equiv\frac{(-1)^{d-1}\tilde{c}_{ab}}{f}
\biggl\{\sum_{k=0}^d(-1)^k A B\frac{\det(e_{I\setminus\{i_k\}})x_{i_k}f_{i_k}
\hat{x}_I\Omega}{f_{i_0}\cdots f_{i_{d}}}\biggr\}_I,$$
where $I=(i_0,\dots,i_d)$. But  the identity 
$\sum_{k=0}^d(-1)^k\det(e_{I\setminus\{i_k\}})e_{i_k}=0$ holds
 and gives  an Euler formula 
$c_I^\beta f=\sum_{k=0}^d(-1)^k\det(e_{I\setminus\{i_k\}})x_{i_k}f_{i_k}$ \cite{bc}.
Thus, 
$$\delta(\psi)\equiv c_{ab}
\fracwithdelims\{\}{A B c_I^\beta
\hat{x}_I\Omega}{f_{i_0}\cdots f_{i_d}}_I,$$
where $c_{ab}=(-1)^{d-1}\tilde{c}_{ab}$.  
\end{pf}

As in the classic case \cite{ps}, we will go further to relate the multiplicative structure
on some quotient of the homogeneous ring $S$ to the cup product on 
the middle cohomology of the quasismooth hypersurface $X$ given by a homogeneous polynomial
$f\in S_\beta$. 

\begin{defn}\cite{bc}   
For  $f\in S_\beta$  the {\it Jacobian ideal} 
$J(f)\subset S$ is the ideal 
generated by the partial derivatives $\partial f/\partial x_1,\dots, \partial f/\partial x_n$.
Also, the {\it Jacobian ring} $R(f)$ is the quotient ring $S/J(f)$ graded by the Chow group 
$A_{d-1}(\pp)$.
\end{defn}

To show a relation between the cup product and multiplication in  $R(f)$ we will need  two
lemmas. We have the natural map
$S_{(a+1)\beta-\beta_0}\rightarrow H^{d-1-a,a}(X)$  
that sends $A$ to the corresponding component of $\res(\omega_A)$. 
The map 
$$\res(\omega_{\_})^{d-1-a,a}:R(f)_{(a+1)\beta-\beta_0}\rightarrow H^{d-1-a,a}(X)$$
induced by the above one is well
defined because of the following statement.

\begin{lem}
If $A\in J(f)_{(a+1)\beta-\beta_0}$, then $\res(\omega_A)^{d-1-a,a}=0$.
\end{lem}

\begin{pf}
In case $a=0$ the statement is trivial because $J(f)_{\beta-\beta_0}=0$. Assume that $a>0$.  
By Theorem \ref{t:res}, since $A\in J(f)$, it suffices to show that 
$\fracwithdelims\{\}{f_jK_{i_a}\cdots K_{i_0}\Omega}{f_{i_0}\cdots f_{i_a}}_{i_0\dots i_a}$
in $C^a({\cal U}|_X,\Omega_X^{d-1-a})$
is a \v{C}ech coboundary for one of the partial derivatives $f_j=\partial f/\partial x_j$. 
We have
\begin{align*}
K_jK_{i_{a}}\cdots K_{i_0}(\dd f\wedge\Omega)=&(-1)^{a+2}\dd f\wedge K_jK_{i_{a}}\cdots
K_{i_0}\Omega+(-1)^{a+1}f_jK_{i_{a}}\cdots K_{i_0}\Omega \\
&+\sum^a_{k=0}(-1)^k f_{i_k}
K_jK_{i_{a}}\cdots\widehat{K_{i_k}}\cdots K_{i_0}\Omega.
\end{align*}
But  $\dd f\wedge\Omega\equiv0$ mod $f$ by equation (\ref{e:equ}), and 
$\dd f=0$ on the hypersurface $X$. Therefore, on $X$ we have the identity:
$$f_jK_{i_a}\cdots K_{i_0}\Omega=\sum_{k=0}^a (-1)^{a+k}
f_{i_k}K_jK_{i_{a}}\cdots\widehat{K_{i_k}}\cdots K_{i_0}\Omega.$$
Hence,
$\fracwithdelims\{\}{f_jK_{i_a}\cdots K_{i_0}\Omega}{f_{i_0}\cdots f_{i_a}}_{i_0\dots i_a}$
is the image of
$(-1)^a\fracwithdelims\{\}{K_j K_{j_{a-1}}\cdots K_{j_0}\Omega}{f_{j_0}\cdots f_{j_{a-1}}}
_{j_0\dots j_{a-1}}$
under the \v{C}ech coboundary map
$C^{a-1}({\cal U}|_X,\Omega^{d-1-a}_X)\rightarrow C^{a}({\cal U}|_X,\Omega^{d-1-a}_X)$.
\end{pf}

Consider the map $S_{(d+1)\beta-2\beta_0}\rightarrow H^{d,d}(\pp)$  that sends
a polynomial $h$ to the class in $H^{d,d}(\pp)$  represented
by the cocycle 
$$\fracwithdelims\{\}{h c_I^\beta\hat{x}_I\Omega}{f_{i_0}\cdots f_{i_d}}
_I\in\check{H}^d({\cal U},\Omega_{\pp}^d)$$ 
as in Theorem \ref{t:cup}.
This induces the map $\lambda:R(f)_{(d+1)\beta-2\beta_0}\rightarrow H^{d,d}(\pp)$
well defined by the following statement.

\begin{lem}\label{l:l2}
If $h\in J(f)$, then $\fracwithdelims\{\}{h c_I^\beta\hat{x}_I\Omega}{f_{i_0}\cdots f_{i_d}}_I$
is a \v{C}ech coboundary.
\end{lem}

\begin{pf}
We can assume that $h$ is one of the partial derivatives $f_j=\partial f/\partial x_j$.
Let $I$ be the ordered subset $\{i_0,\dots,i_d\}\subset\{1,\dots,n\}$. Then 
the equality 
$$c_I^\beta e_j+\sum_{k=0}^d(-1)^{k+1}c_{\{j\}\cup I\setminus\{i_k\}}^\beta e_{i_k}=0$$
(here,  $\{j\}\cup I\setminus\{i_k\}$ is the ordered set $\{j,i_0,\dots,\widehat{i_k},\dots,i_d\}$)
holds and gives the Euler formula \cite{bc}
$$(c_I^\beta b_j+\sum_{k=0}^d(-1)^{k+1}c_{\{j\}\cup I\setminus\{i_k\}}^\beta b_{i_k})f=
c_I^\beta x_j f_j+\sum_{k=0}^d(-1)^{k+1}c_{\{j\}\cup I\setminus\{i_k\}}^\beta x_{i_k}f_{i_k},$$
where the numbers $b_i$ are determined by $\beta=[\sum_{i=1}^n b_i D_i]$. But the number
$c_I^\beta b_j+\sum_{k=0}^d(-1)^{k+1}c_{\{j\}\cup I\setminus\{i_k\}}^\beta b_{i_k}$ 
is the determinant of a matrix with the same two rows $(b_j,b_{i_0},\dots,b_{i_d})$, 
so it vanishes. Using the above Euler formula, we see that  under the \v{C}ech coboundary map
$C^{d-1}({\cal U},\Omega^{d}_\pp)\rightarrow C^{d}({\cal U},\Omega^{d}_\pp)$, the cocycle
$$\fracwithdelims\{\}{f_j c_I^\beta\hat{x}_I\Omega}{f_{i_0}\cdots f_{i_d}}_I=
\fracwithdelims\{\}{\sum_{k=0}^d (-1)^k c_{\{j\}\cup I\setminus\{i_k\}}^\beta f_{i_k}
\hat{x}_{\{j\}\cup(I\setminus\{i_k\})}\Omega}{f_{i_0}\cdots f_{i_d}}_I$$
is the image of 
$\fracwithdelims\{\}{c_{\{j\}\cup J}^\beta\hat{x}_{\{j\}\cup J}\Omega}{f_{j_0}\cdots 
f_{j_{d-1}}}_J$, 
where $J=\{j_0,\dots,j_{d-1}\}$ and $\{j\}\cup J$ is the ordered set $\{j,j_0,\dots,j_{d-1}\}$.
\end{pf}

As a consequence of  Theorem \ref{t:cup} and the above two lemmas, we have proved

\begin{thm}\label{t:dia}
Let $X\subset\pp$ be a quasismooth hypersurface 
defined by $f\in S_\beta$, and suppose $a+b=d-1$. 
Then the diagram
$$\minCDarrowwidth0cm
\newcommand{\topar}{@>c_{ab}\cdot{\rm multiplication}>>}
\newcommand{\botar}{@>\cup>> H^{d-1,d-1}(X) @>\delta>>}
\begin{CD}
R(f)_{(a+1)\beta-\beta_0}\times R(f)_{(b+1)\beta-\beta_0} @>c_{ab}\cdot{\rm multiplication}>>
R(f)_{(d+1)\beta-2\beta_0} \\
@V{\scriptstyle\res(\omega_{\_})^{ba}\times\res(\omega_{\_})^{ab}}VV 
  @V\lambda VV    \\
H^{b,a}(X)\times H^{a,b}(X)\;\:\stackrel{\cup}{\longrightarrow}\;\: H^{d-1,d-1}(X) \;\:
\stackrel{\delta}{\longrightarrow}
\negthickspace\negthickspace\negthickspace\negthickspace\negthickspace\negthickspace
\negthickspace\negthickspace\negthickspace\negthickspace\negthickspace\negthickspace
\negthickspace\negthickspace\negthickspace\negthickspace\negthickspace\negthickspace
\negthickspace\negthickspace\negthickspace\negthickspace\negthickspace\negthickspace
\negthickspace\negthickspace\negthickspace\negthickspace\negthickspace\negthickspace
\negthickspace\negthickspace\negthickspace\negthickspace\negthickspace\negthickspace
\negthickspace\negthickspace\negthickspace\negthickspace\!
@.  H^{d,d}(\pp)
\end{CD}
$$
commutes, where $\lambda$ is as defined above and $\delta$ is the Gysin map.
\end{thm}

\section{Cohomology of  regular hypersurfaces}

In this section we will present
an application of the Gysin spectral sequence to computing 
cohomology of regular semiample hypersurfaces in a complete simplicial toric variety
$\pp$. We will obtain an explicit  description of the cup product on some part of 
the middle cohomology of such hypersurfaces. 
Section 3 studied the relation between multiplication in $R(f)$ and the cup product, whereas 
this section will study such a relation of a smaller ring $R_1(f)$ and the cup product.
The rings $R(f)$ and $R_1(f)$ were previously used in \cite{bc} for studying 
the cohomology of ample hypersurfaces.

Let $D=\pp\setminus{\bf T}=\bigcup^{n}_{i=1}D_i$ and let $X\subset\pp$ be a regular
 hypersurface. Then $(X, X\cap D)$ is a toroidal pair \cite[sect. 15]{d} 
and also $X\cap D$ consists of quasismooth components that intersect quasi-transversally.
Therefore, by the results from \cite[sect. 15]{d}, we have
$$\gr_k^W\Omega_X^\hidot(\log(X\cap D))\cong
\bigoplus\begin{Sb}\dim\sigma=k\end{Sb}\Omega^{\hidot-k}_{X\cap V(\sigma)},$$
and the ({\it Gysin}) spectral sequence of this filtered complex \cite[sect. 3.2]{de} 
$$E_1^{pq}=\!{\Bbb H}^{p+q}(X,\gr_{-p}^W\Omega_X^\hidot(\log(X\cap D))\cong\!\!\!
\bigoplus\begin{Sb}\dim\sigma=-p\end{Sb}\!\!\!\!\!\!H^{2p+q}(X\cap V(\sigma))\Rightarrow H^{p+q}
(X\setminus(X\cap D))$$ 
degenerates at $E_2$ and converges to the weight filtration $W_\lodot$ on 
$H^{p+q}(X\setminus(X\cap D))$: 
$$E_2^{pq}=\gr^W_q H^{p+q}(X\setminus(X\cap D)).$$
In particular, note that
\begin{equation}\label{e:H1}
H^1(X)\cong \gr^W_1 H^1(X\cap \ttt).
\end{equation}

Now assume that $X\subset\pp$ is a regular semiample hypersurface.
In this case $X\setminus(X\cap D)=X\cap{\bf T}$ is a nondegenerate affine 
hypersurface in $\bf T$.
Hence, by Lemma \ref{l:lef},
$E_2^{pq}=\gr^W_{q}H^{p+q}(X\cap{\bf T})$ vanishes unless $p+q=d-1$ and $q\geq d-1$, or
$p+q<d-1$ and $q=-2p$.
Therefore, from the Gysin spectral sequence we obtain the following exact sequences. First, for $s$ odd,  
$s<d-1$, we have
$$
0\rightarrow\!\!\!\bigoplus_{\dim\sigma=\frac{s-1}{2}}\!\!\!H^1(X\cap V(\sigma))
\rightarrow\dots\rightarrow\!\!\bigoplus\begin{Sb}\dim\sigma=k\end{Sb}\!\!H^{s-2k}
(X\cap V(\sigma))\rightarrow\dots\rightarrow H^s(X)\rightarrow0,
$$
where the maps are alternating sums of the Gysin morphisms. Next, for $s$ even,  $s<d-1$, we get
$$
0\rightarrow\gr^W_{s}H^{s/2}(X\cap{\bf T})\rightarrow\bigoplus_{\dim\sigma=\frac{s}{2}}
H^0(X\cap V(\sigma))\rightarrow\dots\rightarrow H^s(X)\rightarrow0.
$$
Finally, for $s=d-1$,  
\begin{equation}\label{e:ex}
\dots\rightarrow\bigoplus\begin{Sb}\dim\sigma=k\end{Sb}H^{d-1-2k}(X\cap V(\sigma))
\rightarrow\dots\rightarrow H^{d-1}(X)\rightarrow\gr_{d-1}^W H^{d-1}(X\cap{\bf T})\rightarrow0.
\end{equation}
Similar sequences exist for $s>d-1$ that are  exact except for one term.
We will be mainly concerned with the last exact sequence, which determines the middle cohomology
group of $X$. 

The following fact, contained in
Proposition 5.3 \cite{c3}, characterizes regular hypersurfaces.

\begin{lem}\cite{c3}\label{l:reg}
Let $X\subset\pp$ be a  hypersurface defined by a homogeneous polynomial
$f$. Then $X$ is regular if and only if $x_i\frac{\partial f}{\partial x_i}$,
 $i=1,\dots,n$, do not vanish simultaneously on $\pp$. In this case we call $f$ nondegenerate.
\end{lem}

This lemma shows that in complete simplicial toric varieties regular hypersurfaces are 
quasismooth. We shall  prove a stronger analog of Theorem
\ref{t:dia} for regular semiample hypersurfaces.

In the case $f$ is nondegenerate, the open sets  $\widetilde{U}_i=\{x\in\pp:x_i f_i(x)\ne0\}$ 
cover the toric variety $\pp$. In particular, the open cover
$\widetilde{\cal U}=\{\widetilde{U}_i\}_{i=1}^n$ 
is a refinement of ${\cal U}$ defined in the previous section.

\begin{defn}\cite{bc}\label{d:r1}
Given  $f\in S_\beta$, we get
the ideal quotient 
$$J_1(f)=\langle x_1\partial f/\partial x_1,\dots,x_n
\partial f/\partial x_n\rangle:x_1\cdots x_n$$
(see \cite[p.~193]{clo}) and the ring $R_1(f)=S/J_1(f)$ graded by the Chow group $A_{d-1}(\pp)$.
\end{defn}

To show the relation between multiplication in $R_1(f)$ and the cup product on the hypersurface
defined by $f$ we need some results similar to those in the previous section.

\begin{lem}\label{l:l1r} Let  $f\in S_\beta$ be nondegenerate and let $h\in J_1(f)$, then the cocycle 
$\fracwithdelims\{\}{h c_I^\beta\hat{x}_I\Omega}{f_{i_0}\cdots f_{i_d}}_I$ vanishes in
$\check{H}^d(\widetilde{\cal U},\Omega_{\pp}^d)$. 
\end{lem}

\begin{pf}
The proof of this is similar to the proof of Lemma \ref{l:l2}.
\end{pf}

\begin{thm}\label{t:l2r}
 Let $X\subset\pp$ be a regular semiample hypersurface defined by 
$f\in S_\beta$, and suppose $a+b=d-1$.

{\rm(i)} If $A\in J_1(f)_{(a+1)\beta-\beta_0}$, then $\res(\omega_A)^{b,a}=0$.

{\rm(ii)} 
The map $\res(\omega_{\_})^{b,a}:R_1(f)_{(a+1)\beta-\beta_0}\rightarrow H^{b,a}(X)$ 
is injective, and the natural composition 
$$R_1(f)_{(a+1)\beta-\beta_0}\rightarrow H^{b,a}(X)\rightarrow H^{b,a}(H^{d-1}(X\cap\ttt))$$
is an isomorphism,
 so that we have a natural imbedding 
$\gr_F W_{d-1}H^{d-1}(X\cap\ttt)\hookrightarrow \gr_F H^{d-1}(X)$.
Moreover, we have an isomorphism 
$$H^{b,a}(X)\cong R_1(f)_{(a+1)\beta-\beta_0}\bigoplus\biggl(\sum_{i=1}^n{\varphi_i}_!H^{b-1,a-1}(X\cap 
D_i)\biggr),$$
where ${\varphi_i}_!$ are the Gysin maps for $\varphi_i:X\cap D_i\hookrightarrow X$, and  
$$\res(\omega_{A})^{b,a}\cup {\varphi_i}_!H^{b-1,a-1}(X\cap D_i)=0\mbox{ for all }
A\in R_1(f)_{(a+1)\beta-\beta_0}.$$
\end{thm}

\begin{pf} {\rm(i)}
We will prove the statement using the Poincar\'e duality
$$H^{b,a}(X)\otimes H^{a,b}(X)\rightarrow H^{d-1,d-1}(X),$$ 
where $b=d-1-a$. Since the pairing is nondegenerate, it suffices to show for 
$A\in J_1(f)_{(a+1)\beta-\beta_0}$
that the cup product of $\res(\omega_A)^{b,a}$ with all elements in $H^{a,b}(X)$ vanishes.
For this we need to find the elements that 
span the  group $H^{a,b}(X)$.  

Let $X$ be linearly equivalent to a torus-invariant
divisor $\sum_{i=1}^n a_i D_i$ with $a_i\ge0$, and $\Delta$ be the corresponding polytope
defined by the inequalities $\langle m,e_i\rangle\ge-a_i$. As in \cite{b1}, $S_\Delta$
denotes the subring of ${\Bbb C}[t_0,t_1^{\pm1},\dots,t_d^{\pm1}]$ spanned over $\Bbb C$ by all
monomials of the form $t_0^k t^m=t_0^kt_1^{m_1}\cdots t_1^{m_d}$ where $k\ge0$ and 
$m\in k\Delta$.
We have a natural isomorphism of graded rings (see the proof of Theorem 11.5 \cite{bc})
$$S_\Delta\cong\bigoplus_{k=0}^\infty S_{k\beta}\subset S,$$
sending $t_0^k t^m$ to $\prod^n_{i=1}x_i^{k a_i+\langle m,e_i\rangle}$. 
This  isomorphism induces the bijection
$$S_{k\beta-\beta_0}@>\prod_{i=1}^n x_i>>\langle x_1\cdots x_n\rangle_{k\beta}\cong 
(I_\Delta^{(1)})_k,$$
where $I_\Delta^{(1)}\subset S_\Delta$ is the ideal spanned by all monomials $t_0^k t^m$ such that
$m$ is in the interior of $k\Delta$.

As a consequence of the exact sequence (\ref{e:ex}) and Theorems 7.13, 8.2 \cite{b1}
 we have the  diagram 
\begin{equation}\minCDarrowwidth0.7cm
\begin{CD}
\label{e:dia}
\bigoplus_{i=1}^n H^{a-1,b-1}(X\cap D_i)
@>\oplus {\varphi_i}_!>> H^{a,b}(X)@>>>H^{a,b}(H^{d-1}(X\cap{{\bf T}}))@>>>0 \\
@. @A\res(\omega_{\_})^{ab}AA   @A\res(\tilde{\omega}_{\_})^{ab}AA @.\\
@. S_{(b+1)\beta-\beta_0}@>>> (I_\Delta^{(1)})_{b+1}, @.
\end{CD}
\end{equation}
where  the top row is exact, the right vertical map  is defined by
$$\tilde{\omega}_{t_0^{b+1} t^m}=\frac{t^m}{\tilde{f}(t)^{b+1}}\frac{\dd t_1}{t_1}\wedge
\dots\wedge \frac{\dd t_d}{t_d}$$ 
(here, $\tilde{f}(t)$ is the Laurent polynomial defining the affine hypersurface 
$X\cap {\bf T}$, so that
$t_0\tilde{f}(t)$ corresponds to $f(x)$ under the isomorphism $(S_\Delta)_1\cong S_\beta$)
and $\res^{ab}$   induced by the Poincar\'e residue mapping \cite[sect.~5]{b1}:
$$\res:H^d({\bf T}\setminus(X\cap {\bf T}))\rightarrow H^{d-1}(X\cap {\bf T}).$$
The diagram commutes because the restriction of the form 
$\omega_B=B\Omega/f^{b+1}$, with $B=\prod^n_{i=1}x_i^{(b+1)a_i-1+\langle m,e_i\rangle}$,
 to the torus ${\bf T}$ coincides with 
$\frac{t^m}{\tilde{f}(t)^{b+1}}\frac{\dd t_1}{t_1}\wedge\dots\wedge 
\frac{\dd t_d}{t_d}$ (use the coordinates 
$t_j=\prod_{i=1}^n x_i^{\langle m_j,e_i\rangle}$ on the torus with the fixed integer basis 
$m_1,\dots,m_d$ from Definition \ref{d:om}).  
   
The first row in (\ref{e:dia}) is exact and the composition 
$$S_{(b+1)\beta-\beta_0}@>>>(I_\Delta^{(1)})_{b+1}@>\res(\tilde{\omega}_{\_})^{ab}>>
H^{a,b}(H^{d-1}(X\cap{\bf T}))$$
is surjective by Theorem 8.2 \cite{b1}. Therefore, the group  $H^{a,b}(X)$ is spanned by 
$\res(\omega_B)^{ab}$ for  $B\in  S_{(b+1)\beta-\beta_0}$ and ${\varphi_i}_!(H^{a-1,b-1}(X\cap D_i))$ 
for $i=1,\dots,n$. 

{}From Theorem \ref{t:cup} and Lemma \ref{l:l1r} it follows that
$\res(\omega_A)^{ba}\cup\res(\omega_B)^{ab}=0$ for $A\in J_1(f)_{(a+1)\beta-\beta_0}$ and
 all $B\in S_{(b+1)\beta-\beta_0}$. Also, for any $A\in S_{(a+1)\beta-\beta_0}$ and
$h\in H^{a-1,b-1}(X\cap D_i)$ we have 
\begin{equation}\label{e:van}
\res(\omega_A)^{ba}\cup{\varphi_i}_!h={\varphi_i}_!(\varphi_i^*\res(\omega_A)^{ba}\cup h)
\end{equation}
by the projection formula for Gysin homomorphisms. 
However, $\varphi_i^*\res(\omega_A)^{ba}$ is represented by the restriction of the 
\v{C}ech cocycle
$c_a\fracwithdelims\{\}{AK_{i_a}\cdots K_{i_0}\Omega}{f_{i_0}\cdots f_{i_a}}
_{i_0\dots i_a}\in\check{H}^a(\widetilde{\cal U}|_X,\Omega_X^b)$ 
from Theorem \ref{t:res} to $X\cap D_i$. This restriction vanishes, because, if 
$i\in\{i_0,\dots,i_d\}$ then $\widetilde{U}_{i_0}\cap\dots\cap\widetilde{U}_{i_a}\cap D_i$ is
empty, and, if $i\notin\{i_0,\dots,i_d\}$ then each term in the form $K_{i_a}\cdots K_{i_0}\Omega$
contains $x_i$ or $\dd x_i$. 
Thus, we have shown that the cup product of $\res(\omega_A)^{ba}$,
 $A\in J_1(f)_{(a+1)\beta-\beta_0}$, with all elements in 
$H^{a,b}(X)$ vanishes, and the result follows. 

{\rm(ii)} {}From the diagram (\ref{e:dia}) and part {\rm(i)}, we get a natural map
$R_1(f)_{(a+1)\beta-\beta_0}\rightarrow H^{d-1-a,a}(H^{d-1}(X\cap\ttt))$.
The fact that this map is an isomorphism follows from the proof of Theorem 11.8 \cite{bc}.
Using the diagram (\ref{e:dia}), we can now see that the map 
$\res(\omega_{\_})^{d-1-a,a}:R_1(f)_{(a+1)\beta-\beta_0}\rightarrow H^{d-1-a,a}(X)$ is injective,
and we get the desired description of the middle cohomology group $H^{d-1}(X)$.
By equation (\ref{e:van}), 
$\res(\omega_{A})^{d-1-a,a}\cup {\varphi_i}_!H^{d-2-a,a-1}(X\cap D_i)=0$. 
\end{pf}

Combining Theorem \ref{t:dia} with Lemma \ref{l:l1r} and Theorem \ref{t:l2r}{\rm(i)} we get the 
following result.

\begin{thm}\label{t:dia1}
Let $X\subset\pp$ be  a regular semiample hypersurface
defined by $f\in S_\beta$, and suppose $a+b=d-1$. 
Then the diagram
$$\minCDarrowwidth0cm
\newcommand{\topar}{@>c_{ab}\cdot{\rm multiplication}>>}
\newcommand{\botar}{@>\cup>> H^{d-1,d-1}(X) @>\delta>>}
\begin{CD}
R_1(f)_{(a+1)\beta-\beta_0}\times R_1(f)_{(b+1)\beta-\beta_0} @>c_{ab}\cdot{\rm multiplication}>>
R_1(f)_{(d+1)\beta-2\beta_0} \\
@V{\scriptstyle\res(\omega_{\_})^{ba}\times\res(\omega_{\_})^{ab}}VV 
  @V\lambda VV    \\
H^{b,a}(X)\times H^{a,b}(X)\;\:\stackrel{\cup}{\longrightarrow}\;\: H^{d-1,d-1}(X) \;\:
\stackrel{\delta}{\longrightarrow}
\negthickspace\negthickspace\negthickspace\negthickspace\negthickspace\negthickspace
\negthickspace\negthickspace\negthickspace\negthickspace\negthickspace\negthickspace
\negthickspace\negthickspace\negthickspace\negthickspace\negthickspace\negthickspace
\negthickspace\negthickspace\negthickspace\negthickspace\negthickspace\negthickspace
\negthickspace\negthickspace\negthickspace\negthickspace\negthickspace\negthickspace
\negthickspace\negthickspace\negthickspace\negthickspace\negthickspace\negthickspace
\negthickspace\negthickspace\negthickspace\negthickspace\!
@.  H^{d,d}(\pp)
\end{CD}
$$
commutes, $c_{ab}=\frac{(-1)^{a(a+1)/2+b(b+1)/2+a^2+d-1}}{a!b!}$.
\end{thm}
 
We will finish this section with an explicit procedure of computing  
$$\int_X \res(\omega_A)^{ba}\cup\res(\omega_B)^{ab}.$$ 
To have this we need generalizations of some results in \cite{c3}.

\begin{defn}\cite{c3}
Assume $F_0,\dots,F_d\in S_\beta$ do not vanish simultaneously on a complete toric variety $\pp$.
Then the {\it toric residue map} 
$$\res_F: S_\rho/\langle F_0,\dots,F_d\rangle_\rho\rightarrow {\Bbb C},$$ 
$\rho={(d+1)\beta-\beta_0}$,
is given by the formula $\res_F(H)=\tr_\pp([\varphi_F(H)])$, where 
$\tr_\pp:H^d(\pp,\Omega_{\pp}^d)\rightarrow{\Bbb C}$ is the trace map, and  $[\varphi_F(H)]$  is 
the class represented by the $d$-form $\frac{H\Omega}{F_0\cdots F_d}$ in \v{C}ech cohomology with
respect to the open cover $\{x\in\pp:F_i(x)\ne0\}$.
\end{defn}

\begin{pr} If $F_0,\dots,F_d\in S_\beta$, then there is $J_F\in S_{(d+1)\beta-\beta_0}$ such that
$$\sum_{j=0}^d (-1)^j F_j \dd F_0\wedge\dots\wedge\widehat{\dd F_j}\wedge\dots\wedge\dd F_d=J_F\Omega.$$
Furthermore, if $I=\{i_0,\dots,i_d\}\subset\{1,\dots,n\}$ such that $c_I^\beta\ne0$  
(if $\beta\ne0$, there is at least one such $I$),
then $J_F=\det(\partial F_j/\partial x_{i_k})/c_I^\beta \hat{x}_I$.
The polynomial $J_F$ is called the toric Jacobian of  $F_0,\dots,F_d$.
\end{pr}

\begin{pf} This is essentially Proposition 4.1 in \cite{c3}. To show that $J_F$ coincides with 
the toric Jacobian in \cite{c3} use the Euler formula 
$$c_I^\beta g=\sum_{k=0}^d(-1)^k\det(e_{I\setminus\{i_k\}})x_{i_k}\frac{\partial g}{\partial x_{i_k}} 
\mbox{ for } g\in S_\beta$$
from the proof of Theorem \ref{t:cup}.
\end{pf}

\begin{thm}\label{t:jac}
Let $\pp$ be a complete toric variety, and let $\beta\in A_{d-1}(\pp)$ be semiample.
If $F_0,\dots,F_d\in S_\beta$ do not vanish simultaneously on $\pp$, then:

{\rm(i)} The  toric residue map 
$\res_F: S_\rho/\langle F_0,\dots,F_d\rangle_\rho\rightarrow {\Bbb C}$,
$\rho={(d+1)\beta-\beta_0}$, is an isomorphism.  

{\rm(ii)} If $J_F\in S_{(d+1)\beta-\beta_0}$ is the toric Jacobian of  $F_0,\dots,F_d$, then
$$\res_F(J_F)=d!\vol(\Delta)=\deg(F),$$
where $\Delta$ is the polyhedron associated to a torus-invariant divisor in the equivalence class of 
$\beta$ and $F:\pp\rightarrow{\Bbb P}^d$ is the map defined by $F(x)=(F_0(x),\dots,F_d(x))$.
\end{thm}

\begin{pf} This statement was proved for ample $\beta$ in Theorem 5.1 \cite{c3}, but the proof 
can be applied in our case almost without change.
Indeed, consider the map $F=(F_0,\dots,F_d):\pp\rightarrow{\Bbb P}^d$ given by the sections of 
a semiample line bundle ${\cal O}_\pp(D)$. Since $(D^d)>0$ and $F_0,\dots,F_d$ do not 
vanish simultaneously on $\pp$,
it follows that  $F_0,\dots,F_d$ are linearly independent. We can extend  $F_0,\dots,F_d$ to a basis of
$H^0(\pp,{\cal O}_\pp(D))$ which  gives the associated map $\phi:\pp\rightarrow {\Bbb P}^N$, where 
$N=h^0(\pp,{\cal O}_\pp(D))-1$.
Then the map $F$ factors through the map $\phi$ and a projection 
$$p:{\Bbb P}^N\setminus L\rightarrow{\Bbb P}^d,\quad (y_0,\dots,y_N)\mapsto(y_0,\dots,y_d),$$ 
where $L\subset{\Bbb P}^N$ is a projective subspace defined by $y_0=\dots=y_d=0$. 
By Exercise on p.~73 \cite[sect.~3.4]{f1}, the dimension of the image of $\phi$ is $d$. 
Using a dimension argument, one can show that $p^{-1}(y_0,\dots,y_d)\cap{\rm im}(\phi)$ is 
nonempty. Hence,  $F$ is surjective, and, consequently, generically finite.
Propositions 3.1, 3.2 and 3.3 in  \cite{c3} are still valid in the case $\beta$ is semiample,
because the isomorphism $S_\Delta\cong S_{*\beta}$ holds.  
The rest of the arguments in \cite{c3} applies without change.
\end{pf}

\begin{defn}\cite{bc}
Given  $f\in S_\beta$, let $J_0(f)\subset S$ denote  the ideal generated by
$x_i\partial f/\partial x_i$, $1\le i\le n$, and put $R_0(f)=S/J_0(f)$.   
\end{defn}

\begin{lem}\label{l:ci} If $I=\{i_0,\dots,i_d\}\subset\{1,\dots,n\}$ such that $c_I^\beta\ne0$,
then $x_i\partial f/\partial x_i$, $i\in I$, don't vanish simultaneously on $\pp$, and
$J_0(f)=\langle x_{i_0}\partial f/\partial x_{i_0},\dots,x_{i_d}\partial f/\partial x_{i_d}\rangle$.
\end{lem}

\begin{pf} If $c_I^\beta\ne0$, then $e_{i_0},\dots,e_{i_d}$ span $M_{\Bbb R}$.
{}From the Euler formula 
$c_I^\beta f=\sum_{k=0}^d(-1)^k\det(e_{I\setminus\{i_k\}})x_{i_k}\frac{\partial f}{\partial x_{i_k}}$ 
and Proposition 5.3 \cite{c3} the lemma follows.
\end{pf}

We now return to the calculation of 
$\int_X \res(\omega_A)^{ba}\cup\res(\omega_B)^{ab}$, when $X$ is a regular semiample hypersurface.
Let $F_j=x_j\frac{\partial f}{\partial x_j}$, and let $I=\{i_0,\dots,i_d\}\subset\{1,\dots,n\}$ be such
that  $c_I^\beta\ne0$. Then denote 
$J=\det(\frac{\partial F_j}{\partial x_i})_{i,j\in I}/(c_I^\beta)^2\hat{x}_I$. One can show that 
$J$ does not depend on the choice of $I$.
By Lemma \ref{l:ci}, the polynomials $F_i$, $i\in I$, do not vanish simultaneously on $\pp$, and 
determine the toric residue map $\res_{F_I}$.
{}From the definitions of $\lambda$, $\res_{F_I}$, and Proposition A.1 \cite{c3}, 
we obtain a commutative diagram
\begin{equation}\label{e:res}
\begin{CD}
R_1(f)_{(d+1)\beta-2\beta_0}@>\prod_{i=1}^n x_i>>R_0(f)_{(d+1)\beta-\beta_0} \\
@V\lambda VV @V{c_I^\beta\res_{F_I}}VV \\
H^{d,d}(\pp)@>{\fracwithdelims(){-1}{2\pi \sqrt{-1}}}^d\int_\pp>> {\Bbb C},
\end{CD}
\end{equation}
where the arrow on the top is just the  multiplication.
Using  Theorem \ref{t:jac}, we get the following procedure.
For given $A\in R_1(f)_{(a+1)\beta-\beta_0}$ and $B\in R_1(f)_{(b+1)\beta-\beta_0}$
there is a unique constant $c$ such that 
$$A\cdot B x_1\cdots x_n-cJ\in\langle x_1\partial f/\partial x_1,\dots,x_n
\partial f/\partial x_n\rangle.$$
Then 
$$\int_X \res(\omega_A)^{ba}\cup\res(\omega_B)^{ab}=c(-2\pi \sqrt{-1})^d c_{ab}d!\vol(\Delta_D),$$ 
where $D=\sum_{i=1}^n a_i D_i$ such that $[D]=\beta$.

\section{Cup product on regular semiample threefolds}

 In this section we will completely describe 
the middle cohomology and the cup product on it for a regular semiample hypersurface 
$X\subset\ps$, when $\dim \ps=4$.

It follows from (\ref{e:ex}) that the map
 $\bigoplus_{i=1}^n H^1(X\cap D_i)@>\oplus {\varphi_i}_!>> H^3(X)$ is injective.
Hence, by Theorem \ref{t:l2r}, 
\begin{equation}\label{e:iso}
H^{b,a}(X)\cong R_1(f)_{(a+1)\beta-\beta_0}\bigoplus\biggl(\bigoplus_{i=1}^n {\varphi_i}_! 
H^{b-1,a-1}(X\cap D_i)\biggr),
\end{equation}
where $a+b=3$.
We first determine which of the groups $H^{b-1,a-1}(X\cap D_i)$ vanish.

\begin{lem}\label{l:desc}
 Let $X\subset\ps$, $\dim \ps=4$, be a $\Sigma$-regular semiample hypersurface, and
let $\pi:\ps\rightarrow\psx$ be the morphism associated with $X$. Then

{\rm(i)}
$H^1(X\cap D_i)=0$ unless $\rho_i\subset\sigma$ for some $2$-dimensional cone $\sigma\in\Sigma_X$, and
$\rho_i\notin\Sigma_X(1)$ (so $\rho_i\setminus\{0\}$ lies in the relative interior of $\sigma$).

{\rm(ii)} For $\rho_i\subset\sigma$, such that $\sigma\in\Sigma_X(2)$ and $\rho_i\notin\Sigma_X(1)$, 
we have
$$\pi_i^*:H^1(Y\cap V(\sigma))\cong H^1(X\cap D_i),$$
where  $V(\sigma)=\pi(D_i)$ is the orbit closure corresponding to  $\sigma\in\Sigma_X$, $Y:=\pi(X)$,
and $\pi_i:X\cap D_i\rightarrow Y\cap V(\sigma)$ is the map induced by $\pi$.
\end{lem}

\begin{pf} {\rm(i)}
Applying (\ref{e:H1}) to the regular hypersurface 
$X\cap D_i$ in the toric variety  $D_i$, we have
\begin{equation}\label{e:h1}
H^1(X\cap D_i)\cong \gr^W_1 H^1(X\cap \ttt_{\rho_i}).
\end{equation}
If $\rho_i\in\Sigma_X(1)$ then $X\cap \ttt_{\rho_i}$ is a nondegenerate affine hypersurface in
$\ttt_{\rho_i}$ because of (\ref{e:str}). Hence, 
$$\gr^W_1 H^1(X\cap \ttt_{\rho_i})\cong\gr^W_1 H^1(\ttt_{\rho_i})=0.$$
If $\rho_i$ does not lie in  a cone  $\sigma\in\Sigma_X(2)$, then 
$X\cap \ttt_{\rho_i}$ is empty or a disjoint finite union of $({\Bbb C}^*)^2$, by equation (\ref{e:str}).
 In this case $\gr^W_1 H^1(X\cap\ttt_{\rho_i})$ also vanishes, and the part {\rm(i)} follows.

{\rm(ii)} Suppose $\rho_i\notin\Sigma_X(1)$ is contained in a cone $\sigma\in\Sigma_X(2)$, and let 
$\sigma'\in\Sigma(2)$ be the cone such that $\rho_i\subset\sigma'\subset\sigma$.
Then we get a composition 
$$H^1(Y\cap V(\sigma))@>\pi_i^*>>H^1(X\cap D_i)@>\varphi_{i,\sigma'}^*>>H^1(X\cap V(\sigma')),$$
where $\varphi_{i,\sigma'}:X\cap V(\sigma')\hookrightarrow X\cap D_i$ is the inclusion.
To prove part {\rm(ii)} it suffices to show that this composition is an isomorphism and 
all spaces in the composition are of the same dimension. 
Applying (\ref{e:H1}) to the regular hypersurfaces 
$X\cap V(\sigma')$ in $V(\sigma')$ and $Y\cap V(\sigma)$ in $V(\sigma)$, 
 we get a commutative diagram
$$
\begin{matrix}
H^1(Y\cap V(\sigma))& \cong&\gr^W_1 H^1(Y\cap \ttt_{\sigma})\\
\downarrow&&\downarrow  \\
H^1(X\cap V(\sigma'))&\cong&\gr^W_1 H^1(X\cap \ttt_{\sigma'}),
\end{matrix}
$$ 
where  the vertical arrow on the right is induced by the isomorphism 
$\ttt_{\sigma'}\cong\ttt_{\sigma}$. {}From the diagram we see that the natural map  
$H^1(Y\cap V(\sigma))\rightarrow H^1(X\cap V(\sigma'))$  is an isomorphism.
On the other hand, since $X\cap \ttt_{\rho_i}\cong (Y\cap \ttt_\sigma)\times{\Bbb C}^*$ because of 
(\ref{e:str}), it follows from (\ref{e:h1}) that   
$$H^1(X\cap D_i)\cong\gr^W_1 H^1((Y\cap \ttt_\sigma)\times{\Bbb C}^*)\cong\gr^W_1 H^1(Y\cap \ttt_\sigma)$$
by the K\"unneth isomorphism. Thus, the dimensions of spaces $H^1(X\cap D_i)$ and 
$H^1(Y\cap V(\sigma))$
coincide. This finishes the proof of part  {\rm(ii)}.
\end{pf}

The above lemma relates the nonvanishing groups $H^1(X\cap D_i)$ to the middle cohomologies of regular
ample  hypersurfaces in  2-dimensional toric varieties. 
Using (\ref{e:iso}) and Theorem \ref{t:l2r}, we can now give a complete algebraic description of the middle
cohomology group $H^3(X)$.
Let $S(V(\sigma))={\Bbb C}[y_\gamma:\sigma\subset\gamma\in\Sigma_X(3)]$ be  the coordinate ring of 
the 2-dimensional complete toric variety 
$V(\sigma)\subset\psx$, and let $f_\sigma\in S(V(\sigma))_{\beta^\sigma}$ 
denote the polynomial defining the hypersurface $Y\cap V(\sigma)$ in $V(\sigma)$. Then, as in Definition
 \ref{d:r1}, we have the ideal $J_1(f_\sigma)$ in $S(V(\sigma))$ and the quotient ring 
$R_1(f_\sigma)=S(V(\sigma))/J_1(f_\sigma)$. By Theorem \ref{t:l2r}(ii), we have an isomorphism 
$$H^{2-a,a-1}(Y\cap V(\sigma))\cong R_1(f_\sigma)_{a \beta^\sigma-\beta_0^\sigma},$$ 
where  $\beta_0^\sigma=\deg(\prod_\gamma y_\gamma)\in A_1(V(\sigma))$. We can now state our first 
main result of this section.

\begin{thm}\label{t:m1}
Let $X\subset\ps$, $\dim\ps=4$, be  a regular semiample hypersurface
defined by $f\in S_\beta$.  
Then there is a natural isomorphism
$$H^{3-a,a}(X)\cong R_1(f)_{(a+1)\beta-\beta_0}
\bigoplus\biggl(\bigoplus\begin{Sb}\sigma
\in\Sigma_X(2)\end{Sb}(R_1(f_\sigma)_{a\beta^\sigma-\beta_0^\sigma})
^{n(\sigma)}\biggr),$$
where $n(\sigma)$ is the number of cones $\rho_i$ such that $\rho_i\subset\sigma$ and $\rho_i\notin
\Sigma_X(1)$. 
\end{thm}

\begin{rem} As we mentioned in introduction and in Remark \ref{r:mir}, in the Batyrev mirror construction 
\cite{b2} a MPCP-desingularization $\widehat{Z}$ of an ample  Calabi-Yau hypersurface of a toric Fano variety 
$\pp_\Delta$, corresponding to a reflexive polytope $\Delta$,  is a regular semiample hypersurface.
In Corollary 4.5.1 \cite{b2} Batyrev calculated the Hodge number 
$$h^{2,1}(\widehat{Z})=l(\Delta)-5-\sum_{{\rm codim}\theta=1}l^*(\theta)+\sum_{{\rm codim}\theta=2}
l^*(\theta)l^*(\theta^*),$$
where $\theta$ is a face of $\Delta$, $\theta^*$ is the corresponding dual face of the dual reflexive polyhedron
$\Delta^*$, and $l(\Gamma)$ (resp.~$l^*(\Gamma)$) denotes the number of integer (resp.~interior integer) points
in $\Gamma$.
We can compare this number with the algebraic description of $H^{2,1}(\widehat{Z})$ in the above theorem.
{}From Theorem \ref{t:l2r} we know that $\dim R_1(f)_{2\beta-\beta_0}=h^{2,1}(\widehat{Z}\cap\ttt)$,
which is equal to  $l(\Delta)-5-\sum_{{\rm codim}\theta=1}l^*(\theta)$ by Theorem 4.3.1 \cite{b2}.
The number $l^*(\theta^*)$ is equal to $n(\sigma)$ of the above theorem for the cone $\sigma$, corresponding
to the face $\theta$ of $\Delta$. 
And finally, one can verify that $\dim R_1(f_\sigma)_{\beta^\sigma-\beta_0^\sigma}$ corresponds to $l^*(\theta)$.
We can now see how the formula for the Hodge number  $h^{2,1}(\widehat{Z})$ is related to our algebraic 
description.
\end{rem}

The next thing we want to do is to compute the cup product on $H^3(X)$ in terms of the algebraic description
in the above theorem.  
To compute this cup product we need one topological result.

\begin{lem}\label{l:com}
Let $K$, $L$ be  subvarieties of a compact V-manifold $M$, which intersect quasi-transversally,
and suppose that $K$, $L$ and $K\cap L$ are  compact V-manifolds.
Then the diagram 
$$
\begin{CD}
H^\hidot(K)@>i_!>>H^\hidot(M)\\
 @V{i'}^*VV @V j^*VV \\
H^\hidot(K\cap L)@>\alpha\cdot{j'}_!>> H^\hidot(L), 
\end{CD}
$$
commutes,
where $i, j, i', j'$ are inclusions, and the constant $\alpha$ satisfies   $[K]\cup[L]=\alpha[K\cap L]$
for fundamental cohomology classes of $K$, $L$ and $K\cap L$ in $M$.
\end{lem}

\begin{pf} The arguments are the same as in the proof of Proposition 10.9 \cite[{\rm VIII}]{do}. 
The only difference is caused by the difference between  $[K]\cup[L]$ and $[K\cap L]$ 
(in the smooth case we won't see this difference).
\end{pf}

\begin{ex} A simple nontrivial example of the above lemma occurs when $M$ is a 2-dimensional toric variety and 
$K$, $L$ are irreducible torus-invariant divisors, intersecting in a point.  
In this case we have to compare the composition of maps
$H^0(K)\stackrel{i_!}{\rightarrow}H^2(M)\stackrel{j^*}{\rightarrow}H^2(L)$ with 
$H^0(K)\stackrel{{i'}^*}{\rightarrow}H^0(K\cap L)\stackrel{{j'}_!}{\rightarrow} H^2(L)$. Since $H^0(M)\stackrel{i^*}{\cong}H^0(K)$ and 
$H^2(L)\stackrel{j_!}{\cong}H^4(M)$ are isomorphisms, it suffices to compare $j_!j^*i_!i^*=[K]\cup[L]\cup\_$ with 
$j_!{j'}_!{i'}^*i^*=[K\cap L]\cup\_$ on $H^0(M)$. The difference between $[K]\cup[L]$ and $[K\cap L]$ can be easily
determined by means of the ring isomorphism $A^\hidot(M)\otimes{\Bbb C}\cong H^{2\hidot}(M)$ \cite[sect.~10]{d}, which sends 
a cycle class of a subvariety $V$ to its fundamental cohomology class $[V]$ in $M$. 
\end{ex}

Equation (\ref{e:iso}) provides a  description of the middle cohomology group $H^3(X)$.
We first show where the cup product on $H^3(X)$ vanishes.

\begin{lem}\label{l:cup1}
${\varphi_i}_!H^1(X\cap D_i)\cup{\varphi_j}_!H^1(X\cap D_j)=0$, $i\ne j$, unless $\rho_i$, $\rho_j$ 
span a cone $\sigma'\in\Sigma$ contained in a $2$-dimensional cone of $\Sigma_X$.
\end{lem}

\begin{pf} 
 By the projection formula for Gysin homomorphisms, we know that 
$${\varphi_i}_!\_\cup{\varphi_j}_!\_={\varphi_j}_!(\varphi_j^*{\varphi_i}_!\_\cup\_).$$
By Lemma \ref{l:com}, for $i\ne j$ we have a commutative diagram
\begin{equation}\label{e:gys}
\minCDarrowwidth0.6cm
\begin{CD}
H^1(X\cap D_i)@>{\varphi_i}_!>>H^3(X)\\
 @V\varphi_{ij}^*VV @V\varphi_j^*VV \\
H^1(X\cap D_i\cap D_j)@>\alpha_{ij}{\varphi_{ji}}_!>> H^3(X\cap D_j), 
\end{CD}
\end{equation}
where $\varphi_{ij}:X\cap D_i\cap D_j\hookrightarrow X\cap D_i$ is the inclusion map and $\alpha_{ij}$
is an appropriate constant.
Hence, it suffices to show that $H^1(X\cap D_i\cap D_j)=0$. This is so, if $\rho_i$ and $\rho_j$ 
do not span a cone in $\Sigma$, because $D_i\cap D_j$ is an empty set in this case.
If $\rho_i$ and $\rho_j$  span a cone  $\sigma'\in\Sigma$, then  $D_i\cap D_j=V(\sigma')$. 
Applying (\ref{e:H1}) to the regular hypersurface $X\cap V(\sigma')$ in $V(\sigma')$,
we see 
$$H^1(X\cap V(\sigma'))\cong\gr^W_1 H^1(X\cap \ttt_{\sigma'}).$$
On the other hand, if $\sigma'$ is not contained in a 2-dimensional cone of $\Sigma_X$, then
$X\cap \ttt_{\sigma'}$ is empty or a disjoint finite union of ${\Bbb C}^*$, by equation (\ref{e:str}).
In this case
$\gr^W_1 H^1(X\cap \ttt_{\sigma'})=0$, and the result follows.
\end{pf}

{}From the above result and Lemma \ref{l:desc} we can see that the cup product of two different spaces 
${\varphi_i}_!H^1(X\cap D_i)$ and ${\varphi_j}_!H^1(X\cap D_j)$ vanishes unless we assume
that $\rho_i\setminus\{0\}$ and $\rho_j\setminus\{0\}$ lie in the relative interior of a 2-dimensional cone  
$\sigma\in\Sigma_X$ and $\rho_i$, $\rho_j$ span a cone $\sigma'\in\Sigma$:

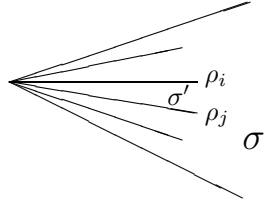
\setlength{\unitlength}{1cm}
\begin{picture}(8,3.3)
\put(2,1.9){\line(3,1){3.2}}
\put(2,1.9){\line(5,1){2.3}}
\put(2,1.9){\line(1,0){2.5}}
\put(4.6,1.9){$\rho_i$}
\put(4.1,1.6){$\sigma'$}
\put(4.6,1.4){$\rho_j$}
\put(5.1,1){\LARGE$\sigma$}
\put(2,1.9){\line(6,-1){2.5}}
\put(2,1.9){\line(3,-1){2.3}}
\put(2,1.9){\line(2,-1){3.1}}
\end{picture}

In this case, by Lemma \ref{l:desc}(ii), we have natural isomorphisms
$${\varphi_i}_!H^1(X\cap D_i)\cong H^1(Y\cap V(\sigma))\cong{\varphi_j}_!H^1(X\cap D_j),$$ 
which provide a natural way to compute the cup product on different spaces:

\begin{lem}\label{l:cup2}
If $\rho_i\ne\rho_j$, not belonging to $\Sigma_X$,  span a 
cone $\sigma'\in\Sigma$  contained in a cone  $\sigma\in\Sigma_X(2)$, then
$${\varphi_i}_!\pi_i^*l_1\cup{\varphi_j}_!\pi_j^*l_2=\frac{{\varphi_{\sigma'}}_!\pi_{\sigma'}^*
(l_1\cup l_2)}{\mult(\sigma')}$$ for $l_1,l_2\in H^1(Y\cap V(\sigma)),$
where $\pi_{\sigma'}:X\cap V(\sigma') \rightarrow Y\cap V(\sigma)$ is the projection and 
$\varphi_{\sigma'}:X\cap V(\sigma')\hookrightarrow X$ is the inclusion.
\end{lem}

\begin{pf}
 Suppose that $\rho_i$ and $\rho_j$  span a cone  $\sigma'\in\Sigma$ contained in 
$\sigma\in\Sigma_X(2)$. Then, using (\ref{e:gys}) and the projection formula, for 
$l_1,l_2\in H^1(Y\cap V(\sigma))$ we compute 
\begin{align*}
&{\varphi_i}_!\pi_i^*l_1\cup{\varphi_j}_!\pi_j^*l_2={\varphi_j}_!(\varphi_j^*{\varphi_i}_!\pi_i^*l_1
\cup\pi_j^*l_2)={\varphi_j}_!(\alpha_{ij}{\varphi_{ji}}_!\varphi_{ij}^*\pi_i^*l_1\cup\pi_j^*l_2)\\
&={\varphi_j}_!{\varphi_{ji}}_!(\alpha_{ij}\varphi_{ij}^*\pi_i^*l_1\cup\varphi_{ji}^*\pi_j^*l_2)
=\alpha_{ij}{\varphi_j}_!{\varphi_{ji}}_!\varphi_{ij}^*\pi_i^*(l_1\cup l_2)={\varphi_j}_!\varphi_j^*
{\varphi_i}_!\pi_i^*(l_1\cup l_2).
\end{align*}
We want to compare the map 
$${\varphi_j}_!\varphi_j^*{\varphi_i}_!\pi_i^*:H^2(Y\cap V(\sigma))\rightarrow H^6(X)$$
 with the map
$${\varphi_{\sigma'}}_!\pi_{\sigma'}^*:H^2(Y\cap V(\sigma))\rightarrow H^6(X).$$
These are the linear maps between 1-dimensional spaces, so they differ by a multiple of a constant.
We will determine this constant using the two commutative diagrams:
$$
\minCDarrowwidth0.6cm
\begin{CD}
H^2(\psx)@>\pi^*>>H^2(\ps)@.@.@.\\
@VVV @VV\varphi_i^*V  @. @. @.\\
H^2(V(\sigma))@>\pi^*_i>>H^2(D_i)@>{\varphi_i}_!>>H^4(\ps)@>\varphi_j^*>>H^4(D_j)@>{\varphi_j}_!>>H^6(\ps)\\
@VVV @VVV @VVV @VVV @VV i^*V \\
H^2(Y\cap V(\sigma))@>\pi^*_i>>H^2(X\cap D_i)@>{\varphi_i}_!>>H^4(X)@>\varphi_j^*>>H^4(X\cap D_j)@>
{\varphi_j}_!>>H^6(X)\\
@. @. @. @. @VV i_!V \\
@. @. @. @. H^8(\ps),
\end{CD}
$$
$$
\minCDarrowwidth0.6cm
\begin{CD}
H^2(\psx)@>\pi^*>>H^2(\ps)@.@. \\
@VVV @VV\varphi_{\sigma'}^*V  @.@. \\
H^2(V(\sigma))@>\pi^*_{\sigma'}>>H^2(V(\sigma'))@>{\varphi_{\sigma'}}_!>>H^6(\ps)@.\\
@VVV @VVV @VV i^*V  \\
H^2(Y\cap V(\sigma))@>\pi^*_{\sigma'}>>H^2(X\cap V(\sigma'))@>{\varphi_{\sigma'}}_!>>H^6(X)@>i_!>>H^8(\ps),
\end{CD}
$$
where the vertical maps are induced by the inclusions. By Lemma \ref{l:com} we had to have some multiplicities
in the above diagrams. These multiplicities are all one because for any $\gamma\in\Sigma$ we have 
$X\cdot V(\gamma)=X\cap V(\gamma)$ in 
$A_\lodot(\ps)$. Indeed, consider a resolution $p:\pp_{\Sigma'}\rightarrow\ps$, corresponding to a nonsingular
subdivision $\Sigma'$ of $\Sigma$. Then, by the proof of Lemma \ref{l:ir}, $p^{-1}(X)\subset\pp_{\Sigma'}$ 
is a regular semiample hypersurface. By the projection formula for cycles, for 
$\gamma'\in\Sigma'(\dim\gamma)$, contained in $\gamma$, we have
$$X\cdot V(\gamma)=p_*(p^*(X)\cdot V(\gamma'))=p_*(p^{-1}(X)\cdot V(\gamma'))=p(p^{-1}(X)\cap V(\gamma'))=
X\cap V(\gamma).$$

We know a nonzero class  $[Y]\in H^2(\psx)$, the fundamental cohomology class of  $Y$ in $\psx$.
Mapping this class to $H^8(\ps)$ in the above two diagrams, we get $[X]\cup[D_i]\cup[D_j]\cup[X]$ and
$[X]\cup[V(\sigma')]\cup[X]$, respectively. Using the  ring isomorphism 
$A^\hidot(\ps)\otimes{\Bbb C}\cong H^{2\hidot}(\ps)$, from Lemma \ref{l:int} we find that 
$[X]\cup[V(\sigma')]\cup[X]$ does not vanish, and, since $D_i\cdot D_j=\frac{1}{\mult(\sigma')}V(\sigma')$
\cite[sect.~5.1]{f1}, it follows that 
$\mult(\sigma'){\varphi_j}_!\varphi_j^*{\varphi_i}_!\pi_i^*={\varphi_{\sigma'}}_!\pi_{\sigma'}^*$ on 
$H^2(Y\cap V(\sigma))$.
\end{pf}

We have computed the cup product of any two different spaces in (\ref{e:iso}). Now we compute the 
cup product on ${\varphi_i}_!H^1(X\cap D_i)$, which does not vanish when 
$\rho_i\setminus\{0\}$ lies in the relative interior of a 2-dimensional cone $\sigma\in\Sigma_X$. In this case 
there are exactly two cones in $\Sigma$, contained in  $\sigma$ and containing $\rho_i$:

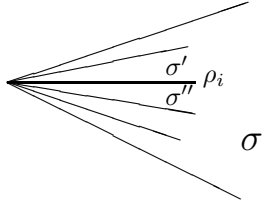
\setlength{\unitlength}{1cm}
\begin{picture}(8,3.3)
\put(2,1.9){\line(3,1){3.2}}
\put(2,1.9){\line(5,1){2.4}}
\put(4.1,2){$\sigma'$}
\put(2,1.9){\line(1,0){2.5}}
\put(4.6,1.9){$\rho_i$}
\put(4.1,1.6){$\sigma''$}
\put(5.1,1){\LARGE$\sigma$}
\put(2,1.9){\line(6,-1){2.5}}
\put(2,1.9){\line(3,-1){2.3}}
\put(2,1.9){\line(2,-1){3.1}}
\end{picture}

In  terms of this, we have

\begin{lem}
 Let $\rho_i\notin\Sigma_X$ be in some $\sigma\in\Sigma_X(2)$ and let $\sigma', \sigma''\in\Sigma(2)$ 
be the two cones, containing $\rho_i$ and contained in $\sigma$. Then
$${\varphi_i}_!\pi_i^*l_1\cup{\varphi_i}_!\pi_i^*l_2=-\frac{\mult(\sigma'+\sigma'')}{\mult(\sigma')\mult(
\sigma'')}{\varphi_{\sigma'}}_!\pi_{\sigma'}^*(l_1\cup l_2)$$ for $l_1,l_2\in H^1(Y\cap V(\sigma))$.
\end{lem}

\begin{pf}
 By the projection formula, we have
$${\varphi_i}_!\pi_i^*\_\cup{\varphi_i}_!\pi_i^*\_={\varphi_i}_!(\varphi_i^*{\varphi_i}_!\pi_i^*\_
\cup\pi_i^*\_)={\varphi_i}_!\varphi_i^*{\varphi_i}_!\pi_i^*(\_\cup\_).$$
As in the proof of the previous lemma,  
we  compare the maps 
${\varphi_i}_!\varphi_i^*{\varphi_i}_!\pi_i^*$ and 
${\varphi_{\sigma'}}_!\pi_{\sigma'}^*$.
Using the arguments of Lemma \ref{l:cup2}, we get
\begin{equation}\label{e:comp}
(X^2\cdot V(\sigma')){\varphi_i}_!\varphi_i^*{\varphi_i}_!\pi_i^*=
(X^2\cdot D_i^2){\varphi_{\sigma'}}_!\pi_{\sigma'}^*
\end{equation}
on $H^2(Y\cap V(\sigma))$.
All we need is to compute the intersection number $(X^2\cdot D_i^2)$.
Take any $m\in M$, such that $\langle m, e_i\rangle\ne0$. The Weil divisor 
$\sum_{j=1}^n \langle m, e_j\rangle D_j$ is equivalent to 0, whence
$$(X^2\cdot D_i^2)=\frac{1}{\langle m, e_i\rangle}(X^2\cdot D_i\cdot(\sum_{j\ne i}-\langle m, e_j\rangle D_j)).
$$
However, $D_i\cdot D_j=\frac{1}{\mult(\gamma)}V(\gamma)$, if $\rho_i$ and $\rho_j$ span a cone 
$\gamma\in\Sigma$,
or $D_i\cdot D_j=0$ otherwise. On the other hand, by Lemma \ref{l:int},  $(X^2\cdot V(\gamma))=0$ unless 
$\gamma$ is contained in $\sigma$. There are exactly two such cones $\sigma'$ and $\sigma''$, contained
in $\sigma$ and containing $\rho_i$. Suppose that $e'$ and $e''$ are the primitive generators of the 
cones $\sigma'$ and $\sigma''$, not lying in $\rho_i$. Then
$$(X^2\cdot D_i^2)=-\frac{\langle m, e'\rangle}{\langle m, e_i\rangle\mult(\sigma')}(X^2\cdot V(\sigma'))-
\frac{\langle m, e''\rangle}{\langle m, e_i\rangle\mult(\sigma'')}(X^2\cdot V(\sigma'')).$$
Since ${\sigma'}^{\perp}={\sigma''}^{\perp}$, equation (\ref{e:int}) shows that 
$(X^2\cdot V(\sigma'))=(X^2\cdot V(\sigma''))$. Also, from \cite[sect.~8.2]{d} it follows that
$\mult(\sigma'+\sigma'')e_i=\mult(\sigma')e''+\mult(\sigma'')e'$. Therefore, 
$$(X^2\cdot D_i^2)=-\frac{\mult(\sigma'+\sigma'')}{\mult(\sigma')\mult(\sigma'')}(X^2\cdot V(\sigma')),$$
and the result follows from equation (\ref{e:comp}).
\end{pf}

We have finished the calculation of the cup product on $H^3(X)$. To state a theorem in a nice form we need to 
define a couple of maps.  
The map $\eta:R_1(f)\rightarrow{\Bbb C}$ is defined 
as $\int_{\ps}\lambda$ on $R_1(f)_{5\beta-2\beta_0}$ (different by a multiple from the map in
(\ref{e:res})), and 0 in all other degrees. Similarly, replacing $\ps$
with $V(\sigma)$ and $f$ with $f_\sigma$, we have the map 
$\eta_\sigma:R_1(f_\sigma)\rightarrow {\Bbb C}$ 
equal to 0 in all degrees except for $3\beta^\sigma-2\beta_0^\sigma$.
Recall also that Theorem \ref{t:m1} gives isomorphism
$$H^{3-a,a}(X)\cong R_1(f)_{(a+1)\beta-\beta_0}
\bigoplus\biggl(\bigoplus\begin{Sb}\sigma
\in\Sigma_X(2)\end{Sb}(R_1(f_\sigma)_{a\beta^\sigma-\beta_0^\sigma})
^{n(\sigma)}\biggr).$$
The following is the description of the cup product on the middle cohomology of the hypersurface.

\begin{thm}
Let $X\subset\ps$, $\dim\ps=4$, be  a regular semiample hypersurface
defined by $f\in S_\beta$.  
If  $A\in R_1(f)_{(a+1)\beta-\beta_0}$, $B\in R_1(f)_{(b+1)\beta-\beta_0}$ are identified with elements of 
$\gr_F H^{3}(X)$ by means of the isomorphism in Theorem \ref{t:m1}, then
$\int_X A\cup B=c_{ab}\eta(A\cdot B)$, where $c_{ab}=\frac{(-1)^{a(a+1)/2+b(b+1)/2+a^2+3}}{a!b!}$.
If we write 
$$( R_1(f)_{a\beta^\sigma-\beta_0^\sigma})^{n(\sigma)}=
\bigoplus_{\sigma\supset\rho_i\notin \Sigma_X}L_a^{\sigma,i},$$ 
where $L_a^{\sigma,i}=R_1(f_\sigma)_{a\beta^\sigma-\beta_0^\sigma}$  
correspond  to the cones $\rho_i$ lying in a $2$-dimensional cone $\sigma\in\Sigma$, 
then for $l_i\in L_a^{\sigma,i}$, $l_i^{\prime}\in L_b^{\sigma,i}$, $l_{j}\in L_b^{\sigma,j}$ 
(identified with elements of $\gr_F H^{3}(X)$) we have
$$\int_X l_i\cup l_{j}=(-1)^{a-1}\frac{\eta_\sigma(l_i\cdot l_{j})}{\mult(\sigma')}$$
in case $\rho_i$ and $\rho_j$ span a cone $\sigma'\in\Sigma(2)$,
$$\int_X l_i\cup l_{i}^{'}=-\frac{\mult(\sigma'+\sigma'')}{\mult(\sigma')\mult(\sigma'')}(-1)^{a-1}\eta_\sigma
(l_i\cdot l_{i}^{'}),$$
where $\sigma', \sigma''\in\Sigma(2)$ are the two cones, contained in $\sigma$ and containing $\rho_i$.
The cup product vanishes in all other cases.
\end{thm}

\begin{rem}
If $X$ is a MPCP-desingularization $\widehat{Z}$ of an ample  Calabi-Yau hypersurface as in \cite{b2},
then the multiplicity $\mult(\sigma')$ is  1 for all 2-dimensional cones $\sigma'$, by the properties of
a reflexive polytope. Also,  in this case $\mult(\sigma'+\sigma'')=2$ in the above theorem.
\end{rem}

\section{Hodge numbers and a ``counterexample'' in mirror symmetry}

In this section we discuss on what kind of mirror symmetry has to be studied. Mirror symmetry proposes
that if two smooth $m$-dimensional Calabi-Yau varieties $V$ and $V^*$ form a mirror pair, then their Hodge
numbers must satisfy the relations
\begin{equation}\label{e:sym}
h^{p,q}(V)=h^{m-p,q}(V^*), 0\le p,q\le m.
\end{equation}
A construction in \cite{b2}, associated with a pair of reflexive polytopes, satisfies the above equalities
for $q=0,1$ \cite{bd}, even if $V$ and $V^*$ are compact orbifolds (i.e., $V$-manifolds). 
We will compute the Hodge numbers $h^{p,2}$ of a regular semiample hypersurface in a complete
simplicial toric variety $\ps$. 
Then we shall apply our formula to the Batyrev mirror construction \cite{b2},
and check that there is no symmetry for the Hodge numbers of MPCP-desingularizations $\widehat{Z}$ of ample 
Calabi-Yau hypersurfaces $\overline{Z}$ coming from a pair of reflexive polytopes $\Delta$ and $\Delta^*$.
However, Theorem 4.15 \cite{bb} and Theorem 6.9 \cite{bd} show that if these MPCP-desingularizations 
$\widehat{Z}$ are smooth, then the duality (\ref{e:sym}) holds. On the other hand, Theorem 4.15 \cite{bb} 
shows that   (\ref{e:sym}) holds for the string-theoretic Hodge numbers $h_{\rm st}^{p,q}$  of the singular 
ample Calabi-Yau hypersurfaces $\overline{Z}$.
This confirms the idea that mirror symmetry has to be studied for smooth varieties with usual Hodge 
numbers or for singular varieties with string-theoretic Hodge numbers.

In order to compute the Hodge numbers we use the $e^{p,q}$ numbers introduced in \cite{dk}:
$$e^{p,q}(V)=\sum_k (-1)^k h^{p,q}(H^k_c(V)),$$
defined for arbitrary algebraic variety $V$. These numbers satisfy the property 
$e^{p,q}(V)=(-1)^{p+q}h^{p,q}(V)$ if $V$ is a compact orbifold. 

Let $X\subset\ps$, $\dim\ps=d$, be a $\Sigma$-regular semiample hypersurface with the associated map 
$\pi:\ps\rightarrow\psx$, $Y=\pi(X)$, as in Proposition \ref{p:reg}. Using the properties of $e^{p,q}$ 
numbers \cite{dk} and equation (\ref{e:str}), for $p+q>d-1$, $p\ne q$, we compute 
\begin{multline*}
h^{p,q}(X)=(-1)^{p+q}e^{p,q}(X)=(-1)^{p+q}\sum_{\sigma\in\Sigma}e^{p,q}(X\cap\ttt_\sigma) \\
=(-1)^{p+q}\sum_{\begin{Sb}\gamma\in\Sigma_X \\ \gamma\supset\sigma\in\Sigma\end{Sb}}
e^{p-i,q-i}(Y\cap\ttt_\gamma)\cdot e^{i,i}(({\Bbb C}^*)^{\dim\gamma-\dim\sigma}),
\end{multline*}
where the sum is by all $\sigma\in\Sigma$ such that $\gamma\in\Sigma_X$ is the smallest cone containing
$\sigma$. Hence, we get 
$$h^{p,q}(X)=(-1)^{p+q}\sum_{\begin{Sb}\gamma\in\Sigma_X \\
0<k\le\dim\gamma-i\end{Sb}}a_{k}(\gamma)(-1)^{\dim\gamma-k+i}\binom{\dim\gamma-k}{i}
e^{p-i,q-i}(Y\cap\ttt_\gamma),$$
where $a_k(\gamma)$ denotes the number of cones $\sigma\in\Sigma(k)$ such that  $\gamma\in\Sigma_X$
is the smallest cone containing $\sigma$, and $\binom{s}{i}$ is the usual binomial coefficient.
It follows from  a formula in \cite[sect.~3.11]{dk} that in the last sum 
$e^{p-i,q-i}(Y\cap\ttt_\gamma)=0$ unless $(p-i)+(q-i)\le\dim(Y\cap\ttt_\gamma)$ 
(equivalently, $\dim\gamma\le d-1-p-q+2i$). 
We now assume $q=d-3$. Then, for $p>2$, $p\ne d-3$, we have
\begin{align*}
&h^{p,d-3}(X)=(-1)^{p+d-3}\sum_{\begin{Sb}\gamma\in\Sigma_X(l) \\
0<k\le l-i \\ l\le 2-p+2i\end{Sb}}a_{k}(\gamma)(-1)^{l-k+i}
\binom{l-k}{i}e^{p-i,d-3-i}(Y\cap\ttt_\gamma) \\
&=
(-1)^{p+d-3}\sum_{\begin{Sb}\gamma\in\Sigma_X(l) \\
0<k\le l-i\le 2-p+i\\    0\le p-i\le1\end{Sb}}a_{k}(\gamma)(-1)^{l-k+i}
\binom{l-k}{i}e^{d-3-i,p-i}(Y\cap\ttt_\gamma) \\
&=
(-1)^{p+d-3}\biggl(\sum_{\gamma\in\Sigma_X(p)}a_{1}(\gamma)e^{d-2-p,1}(Y\cap\ttt_\gamma)+
\sum_{\gamma\in\Sigma_X(p+1)}a_{1}(\gamma)e^{d-3-p,0}(Y\cap\ttt_\gamma) \\
&+
\sum_{\begin{Sb}\gamma\in\Sigma_X(p+2) \\ 1\le k\le2\end{Sb}}a_{k}(\gamma)(-1)^{k}\binom{p+2-k}{p}
e^{d-3-p,0}(Y\cap\ttt_\gamma)\biggr)
\end{align*}

Let $X$ be linearly equivalent to a torus-invariant divisor $\sum_{i=1}^n b_i D_i$, which gives a polytope
$\Delta$.  By Remark \ref{r:cor}, a cone $\gamma\in\Sigma_X$ corresponds to a face 
$\Gamma_\gamma$ of $\Delta$.
Applying Corollary 5.9 and Proposition 5.8 in \cite{dk}, we get: 
$$e^{d-2-p,1}(Y\cap\ttt_\gamma)=(-1)^{d-p-1}\biggl(l^*(2\Gamma_\gamma)-(d-p+1)l^*(\Gamma_\gamma)-
\sum_{\begin{Sb}\Gamma\subset\Gamma_\gamma \\ {\rm codim}\Gamma=1\end{Sb}}l^*(\Gamma)\biggr)$$
if $\dim\gamma=p$ (here, $l^*(\Gamma)$ is the number of interior integral points in $\Gamma$). Furthermore,
$$e^{d-3-p,0}(Y\cap\ttt_\gamma)=(-1)^{d-p-2}\sum_{\begin{Sb}\Gamma\subset\Gamma_\gamma\\ 
{\rm codim}\Gamma=1\end{Sb}}l^*(\Gamma)\quad\mbox{if }\dim\gamma=p+1,$$
$$e^{d-3-p,0}(Y\cap\ttt_\gamma)=(-1)^{d-p-3}l^*(\Gamma_\gamma)\quad\mbox{if }\dim\gamma=p+2.$$
Substituting these numbers in the above formula, we obtain
\begin{align*}
&h^{p,d-3}(X)=\sum_{\gamma\in\Sigma_X(p)}a_{1}(\gamma)\biggl(l^*(2\Gamma_\gamma)-(d-p+1)l^*(\Gamma_\gamma)-
\sum_{\begin{Sb}\Gamma\subset\Gamma_\gamma \\{\rm codim} \Gamma=1\end{Sb}}l^*(\Gamma)\biggr) \\
&-
\sum_{\gamma\in\Sigma_X(p+1)}a_{1}(\gamma)\biggl(\sum_{\begin{Sb}\Gamma\subset\Gamma_\gamma\\ 
{\rm codim}\Gamma=1\end{Sb}}l^*(\Gamma)\biggr)+\sum_{\gamma\in\Sigma_X(p+2)}(a_{2}(\gamma)-a_{1}(\gamma)(p+1))
l^*(\Gamma_\gamma).
\end{align*}
Simplifying and using Poincar\'e duality, we find  for $p>2$, $p\ne d-3$,
\begin{align*}
&h^{d-1-p,2}(X)=\sum_{\gamma\in\Sigma_X(p)}a_{1}(\gamma)\biggl(l^*(2\Gamma_\gamma)-(d-p+1)
l^*(\Gamma_\gamma)-\sum_{\begin{Sb}\Gamma\subset\Gamma_\gamma \\{\rm codim} \Gamma=1\end{Sb}}l^*(\Gamma)\biggr)\\
&+
\sum_{\gamma\in\Sigma_X(p+2)}l^*(\Gamma_\gamma)\biggl(a_{2}(\gamma)-(p+1)a_{1}(\gamma)-\sum_{\begin{Sb}
\tau\subset\gamma\\ 
{\rm codim}\tau=1\end{Sb}}a_1(\tau)\biggr)
\end{align*}

We now can apply this formula to the mirror construction in \cite{b2}. We recall
that a MPCP-desingularization of a toric variety $\pp_\Delta$, associated with a reflexive polytope $\Delta$, 
is a complete simplicial toric variety, corresponding to a refinement 
$\Sigma$ of the normal fan of $\Delta$ such that the cone generators in the fan $\Sigma$ are exactly $N\cap\Delta^*-\{0\}$
(here, $\Delta^*$ is the dual reflexive polytope).
Notice that if $\widehat{Z}_\Delta$ is a MPCP-desingularization of a regular ample hypersurface 
$\overline{Z}_\Delta\subset\pp_\Delta$, then the toric variety $\pp_\Delta$ coincides  with 
$\pp_{\Sigma_{\hat{Z}_\Delta}}$. We also note that in the above formula the number $a_{1}(\gamma)$, for 
$\gamma\in\Sigma_{\hat{Z}_\Delta}$, is equal to $l^*(\Gamma_\gamma^*)$, where $\Gamma_\gamma^*$ is the dual face of 
$\Delta^*$.  
 
The example we use comes from the reflexive 
polytope $\Delta$ of dimension 7 in $M_{\Bbb R}={\Bbb R}^7$, given by the equations
$$z_i\ge-1, i=1,\dots,7,\quad -2z_1-2z_2-2z_3-2z_4-3z_5-3z_6-3z_7\ge-1.$$
The dual reflexive polytope $\Delta^*$ has vertices at 
$$n_0=(-2,-2,-2,-2,-3,-3,-3), n_1,\dots,n_7,$$
where $n_1,\dots,n_7$ are the standard basis of the lattice $N={\Bbb Z}^7$. Notice that the toric variety $\pp_\Delta$, 
corresponding to the polytope $\Delta$, is the weighted projective space $\pp(1,2,2,2,2,3,3,3)$.
The only integral points in $\Delta^*$ are the vertices and the origin, implying that no subdivision occurs for the 
normal fan of $\Delta$; consequently,  $h^{3,2}(\widehat{Z}_\Delta)=0$ because in the above formula for $h^{3,2}$
all the numbers $a_1$ and $a_2$ vanish. On the other hand, for a dual MPCP-desingularization 
$\widehat{Z}_{\Delta^*}$ of a regular ample hypersurface in $\pp_{\Delta^*}$, the above formula for $h^{d-1-p,2}(X)$ 
with $d=7$ and $p=3$ simplifies to 
$$h^{3,2}(\widehat{Z}_{\Delta^*})=\sum_{\begin{Sb}\Gamma^*\subset\Delta^* \\ \dim\Gamma^*=4\end{Sb}} 
l^*(\Gamma)\cdot l^*(2\Gamma^*),$$
because $l^*(\Gamma^*)=0$ for all faces $\Gamma^*$ of $\Delta^*$.
We want to show that $h^{3,2}(\widehat{Z}_{\Delta^*})$ is positive, which would imply that the duality (\ref{e:sym}) fails 
for the pair $(\widehat{Z}_\Delta,\widehat{Z}_{\Delta^*})$. Indeed, consider the 4-dimensional face $\Gamma^*$ of 
$\Delta^*$ with vertices at  $n_0$, $n_1$, $n_2$, $n_3$ and $n_4$. Then $(-1,-1,-1,-1,-2,-2,-2)$ is the interior integral
point of $2\Gamma^*$. The dual 2-dimensional face $\Gamma$, which has vertices at $(-1,-1,-1,-1,5,-1,-1)$, 
$(-1,-1,-1,-1,-1,5,-1)$, $(-1,-1,-1,-1,-1,-1,5)$,  contains the  integral point $(-1,-1,-1,-1,1,1,1)$ in its relative 
interior. Thus, we have shown that $h^{3,2}(\widehat{Z}_\Delta)\ne h^{3,2}(\widehat{Z}_{\Delta^*})$. This happened 
because the hypersurface $\widehat{Z}_{\Delta^*}$ is singular. 

Other ``counterexamples'' can be easily found in higher dimensions, showing that, in general, 
the duality (\ref{e:sym}) fails for the MPCP-desingularizations of the construction in \cite{b2}. 
We should also point out that there is such ``counterexample''
 for 4-folds (see Example 1.2 \cite{b3}).

\end{document}